\documentclass[12pt,reqno]{article}
\usepackage{graphics,epstopdf}
\usepackage{float}
\usepackage{epsfig}
\usepackage{amssymb,amsmath}
\usepackage{amsfonts}
\usepackage{color}
\usepackage{authblk}
\usepackage{subfigure}
\usepackage{amssymb}
\usepackage{caption}
\usepackage{graphicx}
\usepackage[mathscr]{euscript}
\usepackage{mathrsfs}
\usepackage[utf8]{inputenc}
\usepackage[english]{babel}
\RequirePackage{lineno}

\newtheorem{theorem}{Theorem}[section]
\newtheorem{corollary}{Corollary}[section]

\newtheorem{exmp}{Example}[section]
\newtheorem{counterexmp}{Counterexample}[section]

\newtheorem{definition}[theorem]{Definition}

\catcode`@=11 \@addtoreset{equation}{section}
\renewcommand\theequation{\thesection.\@arabic\c@equation}
\catcode`@=12

\makeatother
\usepackage{hyperref}

\begin{document}
\title{{Ordering Results between Two Extreme Order Statistics} with Heterogeneous Linear Failure Rate Distributed Components}
\author{{\large {CM {\bf Revathi}$^{1}$\thanks {Email address: cmrevathi0401@gmail.com, revathi.cm2022@vitstudent.ac.in},~Rajesh {\bf Moharana}$^{2}$\thanks {Email address (corresponding author): rajeshmoharana31@gmail.com,~rajesh.moharana@vit.ac.in}~and Raju {\bf Bhakta}$^{3}$\thanks {Email address: bhakta.r93@gmail.com,~raju\_bhakta.maths@yahoo.com}}}\\
{\em\small\it $^{1,2}$Department of Mathematics, School of Advanced Sciences, Vellore Institute of Technology, Vellore, Tamil Nadu-632014, India.}\\
{\em\small{\it $^{3}$Department of Mathematics and Basic Sciences, National Institute of Information and Technology University, Neemrana-301705, Rajasthan, India.}}
}
\date{}
\maketitle
\begin{abstract}
Stochastic comparisons of series and parallel systems are important in many areas of engineering, operations research and reliability analysis. These comparisons allow for the evaluation of the performance and reliability of systems under different conditions, and can inform decisions related to system design, probabilities of failure, maintenance and operation. In this paper, we investigate the stochastic comparisons of the series and parallel systems under the assumption that the component lifetimes have independent heterogeneous linear failure rate distributions. The comparisons are established based on the various stochastic orders including magnitude, transform and variability orders. Several numerical examples and counterexamples are constructed to illustrate the theoretical outcomes of this paper. Finally, we summarized our findings with a real-world application and possible future scopes of the present study.
\end{abstract}


\section{Introduction}
In the field of operations research and optimization theory, the formalization and continuous improvement of the methods and techniques are necessary in order to address reliability design problems of many complicated systems. The series and parallel systems are the building blocks of various complex systems. Thus, investigating series and parallel systems are helpful in reliability related optimization theory. For an overview of this problem, one can refer to Kuo and Prasad \cite{kuo2000annotated}. Order statistics, especially the extreme order statistics are useful to describe the lifetimes of series and parallel systems.  Let $X_1,\ldots,X_n$ be $n$ independent random variables, and $X_{1:n}\leq\ldots\leq X_{n:n}$ denotes the corresponding order statistics. The random variable $X_{1:n}=\min \{X_1,\ldots,X_n\}$ is known as the smallest order statistic which represents the lifetime of a series system and the random variable $X_{n:n}=\max \{X_1,\ldots,X_n\}$ is known as the largest order statistic which represents the lifetime of a parallel system. We recall that the series and parallel systems are particular cases of the general $k$-out-of-$n$ system. For more detailed explanation and applications of order statistics, the researcher can refer to Balakrishnan and Rao \cite{balakrishnan19981} and David and Nagaraja \cite{david2004order}.
	
The concepts of stochastic orders play an important role in various areas including operations research, reliability theory and risk analysis etc. For example, let $X$ and $Y$ denote the lifetimes of two systems. Then, if $X$ is smaller than $Y$ in the sense of the usual stochastic order (see Shaked and Shanthikumar \cite{shaked2007stochastic}), then a reliability engineer will prefer to the system with lifetime $Y$. Further, suppose $X$ is smaller than $Y$ with respect to the hazard rate order (see Shaked and Shanthikumar \cite{shaked2007stochastic}), then the system with lifetime $X$ will be preferred than the other. Note that the hazard rate means the instantaneous failure rate of a system or a component. It is worth pointing that the reversed hazard rate order and the likelihood ratio order also play a useful role in various applied areas. In this direction, one may refer to Jewitt \cite{jewitt1991} and Nanda and Shaked \cite{nanda2001hazard}. Due to the importance of the problem of stochastic comparisons of the lifetimes of series and parallel systems, several authors have paid their attention when the components have various statistical distributions. A wide variety of research on the stochastic comparisons between the order statistics are available in the literature, where the component lifetimes follow generalized exponential  distributions (see Balakrishnan et al. \cite{balakrishnan2014stochastic}), Fr{\'e}chet distributions (see Gupta et al. \cite{gupta2015stochastic}),  exponential Weibull distributions (see Fang and Zhang  \cite{fang2015stochastic}),  log-Lindley  distributions (see Chowdhury and Kundu \cite{chowdhury2017stochastic})  etc.  In these references, the comparisons include the usual stochastic order, hazard rate order,  likelihood ratio order, star order, dispersive order, etc. More information on stochastic comparisons can be found in Fang and Zhang  \cite{fang2012new}, Misra and Misra \cite{misra2012new}, Torrado and Kochar \cite{torrado2015stochastic}, Kundu and Chowdhury \cite{kundu2016ordering}, Balakrishnan et al. \cite{balakrishnan2018ordering} and Patra et al. \cite{patra2018some}.   This paper discusses stochastic comparisons of two parallel systems with independent heterogeneous components using the exponentiated Kumaraswamy-G distribution model. It establishes the likelihood ratio order among largest order statistics for heterogeneous multiple-outlier models, and provides numerical examples for illustrations.   Kayal et al. \cite{kayal2022stochastic} presented mainly stochastic comparison of two finite mixture models with respect to usual stochastic order and various examples given to verify the established results.  Kayal et al. \cite{kayal2023some} discussed stochastic comparisons of two parallel systems with independent heterogeneous components using the exponentiated Kumaraswamy-G model,  the likelihood ratio order among largest order statistics for heterogeneous multiple-outlier models, and provided numerical examples for illustrations. Barmalzan \cite{barmalzan2022usual} established some new ordering properties between two parallel systems with random, different samples, in the sense of the usual stochastic and reversed hazard rate orders.  By analyzing the  given stochastic orderings, researchers can decide on the best configuration for their particular requirements.

We now describe the significance of the model considered in this paper. The Linear Failure Rate $\mathcal{(LFR)}$ distribution is a widely used model in reliability engineering to describe components or systems with monotonically linear failure rates. It is particularly useful for modeling systems with constant failure rates over the time, providing a flexible framework for analyzing reliability data and making predictions about system performance.  Some results based on the $\mathcal{LFR}$ distribution  may refer the following literatures.  The $\mathcal{LFR}$ distribution has been used to analyze a lifespan of system or units in reliability or survival analysis. Fundamental structural characteristic of the $\mathcal{LFR}$ distribution with minimum of two independent variables X and Y with exponential ($\alpha$) and Rayleigh ($\beta$) distributions as mentioned below.  A random variable $X$ is said to follow $\mathcal{LFR}$ distribution if its cumulative distribution function $(CDF)$ is given by
\begin{equation}\label{equ1}
F_X(x)=1-\exp\left\{-\left(\alpha x+\frac{\beta}{2}x^2\right)\right\};~x>0,~\alpha,~\beta>0,
\end{equation}
where $\alpha$ and $\beta$ are the shape and scale parameters, respectively.  Here, we use the notation $X \thicksim \mathcal{LFR}(\alpha,\beta)$ if $X$ has the $CDF$ given in (\ref{equ1}). The probability density function $(PDF)$ of the $\mathcal{LFR}$ distribution is
\begin{equation}\label{equ2}
f_X(x)=(\alpha+\beta x)\exp\left\{-\left(\alpha x+\frac{\beta}{2}x^2\right)\right\};~x>0,~\alpha,~\beta>0.
\end{equation}
This distribution has been influenced by its application to human life time data \cite{kodlin1967new}.  This paper discusses the usefulness of a skewed two-parameter distribution in human survival time analysis. It discusses a maximum likelihood technique and damage models that incorporate this distribution. The flexibility of this distribution allows for more accurate representation of real-world data and better understanding of different factors impact on survival times.  Its properties were studied by various authors especially Bain \cite{bain1974analysis}, the study explores life-testing distributions with polynomial hazard rate functions and suggests least squares estimators as a potential parameter estimation method. Results show less bias and efficiency compared to maximum likelihood estimators, and their applicability in complex models is explored. Sen \cite{sen1995inference}, for type II censored samples, the paper investigates the maximum likelihood and least-squares-type estimation of the $\mathcal{LFR}$. It demonstrates the $\mathcal{LFR}$ structural property simplifies EM algorithm application for MLE computations.  An alternative method based on pseudo likelihood maximization is developed for better coverage probabilities in large samples. For instance, suppose, in a telecommunications network, two configurations series and parallel are considered to ensure uninterrupted service. In the series system, communication passes sequentially through three routers $R_1$, $R_2$, and $R_3$, each with heterogeneous linear failure rates of $0.02t$, $0.03t$, and $0.01t$, respectively, where $t$ represents time in hours. The reliability of the series system is the product of the individual reliabilities: $R_s(t)=e^{-0.02t}.e^{-0.03t}.e^{-0.01t}=e^{-0.06t}$, indicating that failure of any single router causes the entire system to fail. Conversely, in a parallel system, communication is split across three redundant paths, each with the same failure rates. The system remains operational as long as at least one path functions, with the failure function expressed as $F_p(t)=(1-e^{-0.02t})(1-e^{-0.03t})(1-e^{-0.01t})$ and the reliability as $R_p(t)=1-F_p(t)$. Stochastically, the series system's reliability declines more rapidly over time due to its sensitivity to individual component failures, whereas the parallel system sustains higher reliability by tolerating individual path failures. This reflects real-world scenarios where redundant paths in communication networks significantly enhance reliability compared to single-path configurations.

In this article, we study the stochastic comparison of series and parallel systems of the $\mathcal{LFR}$ distribution using various stochastic orders. The goal of this paper is to obtain various sufficient conditions, for which two series and parallel systems with the $\mathcal{LFR}$ distributed component lifetimes are comparable in terms of the magnitude, transform, and variability orders. The components are taken as independent and heterogeneous. Suppose there are two aircrafts having four engines, which are connected in either series or parallel. Further, assume that the components' lifetimes follow independent heterogeneous $\mathcal{LFR}$ distributions. Then, the established results of this paper are useful to find the more reliable aircraft.

The paper is organized as follows: In Section $2$, we introduce some basic notations and definitions of stochastic orders that have been used throughout this paper. The results for the usual stochastic ordering,  hazard rate ordering, and likelihood ratio ordering of the series and parallel system for usual stochastic order when components follow independent heterogeneous  $\mathcal{LFR}$ are reported in Section $3$. In Section $4$ stochastic comparison results were studied for dispersive, star and convex orderings. Section $5$ includes various examples and counterexamples to demonstrate the theoretical conclusions. Finally, Section $5$ concludes about the research paper with future work.
	
\section{Preliminaries}
In this section, we present some basic concepts of stochastic orderings. Here, all the random variables are nonnegative and absolutely continuous. The words ``increasing'' and ``decreasing'' respectively mean ``non-decreasing'' and ``non-increasing''. For a continuously differentiable function $\xi(t)$, the first and second order derivatives with respect to $t$ are represented by $\frac{\partial\xi(t)}{\partial t}$ and $\frac{\partial^2\xi(t)}{\partial t^2}$, respectively. Also, we use ``log'' for usual logarithm base ``e''. Further, $\mathbb{R}^+$ denotes the set of all positive real numbers. Let $X$ and $Y$ be two nonnegative independent random variables with the $PDF$s $f_X(\cdot)$ and $g_Y(\cdot)$, the $CDF$s $F_X(\cdot)$ and $G_Y(\cdot)$, the survival functions ($SF$s) $\bar{F}_X(\cdot)\equiv 1-F_X(\cdot)$ and $\bar{G}_Y(\cdot)\equiv 1-G_Y(\cdot)$, the hazard rate functions ($HRF$s) $h_X(\cdot)\equiv f_X(\cdot)/\bar{F}_X(\cdot)$ and $h_Y(\cdot)\equiv g_Y(\cdot)/\bar{G}_{Y}(\cdot)$, respectively.
Various types of stochastic orders have been formed and studied in the literature. The following standard widely recognized definitions might be acquired in Muller \cite{muller2002comparison}, Shaked and Shanthikumar \cite{shaked2007stochastic}, and Li \cite{li2013stochastic}.

\begin{definition}\label{def2.1}
$X$ is stated to be smaller than $Y$ in the sense of the
\begin{enumerate}
\item[(i)] usual stochastic (st) ordering $(X\leq_{st}Y)$, if $\bar{F}_X(x)\leq\bar {G}_Y(x)$, for all $x\in\mathbb{R^+};$

\item[(ii)] hazard rate (hr) ordering $(X\leq_{hr}Y)$, if $\bar{G}_Y(x)/\bar{F}_X(x)$ is increasing, for all $x\in\mathbb{R^+};$ or, equivalently, if $h_X(x)\geq h_Y(x),$ for all $x\in\mathbb{R^+};$
		
		
\item[(iii)] likelihood ratio (lr) ordering $(X\leq_{lr}Y)$, if $g_Y(x)/f_X(x)$ is increasing, for all $x\in\mathbb{R^+}.$
\end{enumerate}
\end{definition}
It is worthwhile to mention here that

\begin{equation*}
X\leq_{lr}Y\Rightarrow X\leq_{hr}Y\Rightarrow X\leq_{st}Y.
\end{equation*}
Next, we recall some transform and variability orders that are used in this article.

\begin{definition}\label{def2.2}
$X$ is stated to be smaller than $Y$ in the sense of
\begin{enumerate}
\item[(i)] the dispersive (disp) order $(X\leq_{disp}Y)$, if ${F_{X}^{-1}(\beta)}- {F_{X}^{-1}(\alpha)}\leq {{G_{Y}^{-1}(\beta)}-{G_{Y}^{-1}(\alpha)}}$ whenever $0<{\alpha}\leq{\beta}<1$;
or, equivalently, if {$G_{Y}^{-1}(p)-F_{X}^{-1}(p)$} is increasing in $p\in(0,1)$;
\item[(ii)] the star (*) order $(X\leq_{*}Y)$, if {$\frac{G_{Y}^{-1}(F_{X}(x))}{x}$} is increasing in $x$, for all $x\in\mathbb{R^+};$ or, equivalently, if {$\frac{G_{Y}^{-1}(p)}{F_{X}^{-1}(p)}$} is increasing in {$p\in(0,1)$};
\item[(iii)] the convex transform (c) order $(X\leq_{c}Y)$, if {${G_{Y}^{-1}(F_{X}(x))}$} is convex in $x$, for all $x\in\mathbb{R^+}$;
\item[(iv)]the Lorenz (Lorenz) order ($X\leq_{Lorenz}Y$), if $\frac{1}{E(X)}\int_{t}^{1}F^{-1}_{X}(x)dx\leq\frac{1}{E(Y)}\int_{t}^{1}G^{-1}_{Y}(x)dx$, for all $t\in[0,1]$ and for which the expectations are exist.
\end{enumerate}
\end{definition}
It is well known that
\begin{equation*}
X\leq_{c}Y\Rightarrow X\leq_{*}Y\Rightarrow X\leq_{Lorenz}Y.
\end{equation*}
For a comprehensive exploration of various stochastic orderings, refer to Shaked and Shanthikumar \cite{shaked2007stochastic} and Kleiber \cite{kleiber2002variability}.
	
\section{Results based on Magnitude Orders}\label{section3}
In this section, we derived various ordering outcomes for the lifetimes of series systems in the sense of the usual stochastic, hazard rate, and likelihood ratio orderings. Let $\left\{X_{1},\ldots,X_{n}\right\}$ and $\left\{Y_{1},\ldots,Y_{n}\right\}$ be two independent and identically distributed random samples that follow the $\mathcal{LFR}$ distribution with specified parameters. The following theorem deals with the usual stochastic ordering of series systems.

\begin{theorem}\label{(3.1)}
Let $\{X_1,\ldots,X_n\}$ and $\{Y_1,\ldots,Y_n\}$ be two sets of random variables where $X_k\thicksim LFR(\alpha_k,\beta_k)$ and $Y_k\thicksim LFR(\alpha^*_k,\beta^*_k)$ for $k=1,\ldots,n,$ respectively. If $\alpha_k\geq\alpha^*_k$ and $\beta_k\geq\beta^*_k$, then $X_{1:n}\leqslant_{st}Y_{1:n}$.
\end{theorem}
{\bf Proof:} The $CDF$s of $X_{1:n}$ and $Y_{1:n}$ are given by
\begin{equation}\label{(eqn.3.1)}
F_{X_{1:n}}(x)=1-\exp\left\{-\sum\limits_{k=1}^n\big(\alpha_k x+\frac{\beta_k}{2}x^2\big)\right\}
\end{equation}
and
\begin{equation}\label{(eqn.3.2)}
G_{Y_{1:n}}(x)=1-\exp\left\{-\sum\limits_{k=1}^n\big(\alpha^*_kx+\frac{\beta^*_k}{2}x^2\big)\right\},
\end{equation}
respectively. Now, using Part $(i)$ of Definition \ref{def2.1}, we need to prove that $\bar{F}_{X_{1:n}}(x) \leq \bar{G}_{Y_{1:n}}(x)$. {In doing so, it is required to show that}
\begin{eqnarray}\label{eq3.3}
\alpha_kx+\frac{\beta_k}{2}x^2\geq\alpha^*_kx+\frac{\beta^*_k}{2}x^2.
\end{eqnarray}
{The above inequality given in (\ref{eq3.3}) is satisfied} under the following conditions
\begin{equation*}
\alpha_k\geq\alpha^*_k~{\mbox{and}}~\beta_k\geq\beta^*_k,
\end{equation*}
{which implies that $X_{1:n}\leqslant_{st}Y_{1:n}$. Hence, the result follows.}
\hfill\(\blacksquare\)

To validate the result in Theorem \ref{(3.1)}, we present Example \ref{example5.1} and Counterexample \ref{counexam5.1} in Section \ref{section5}. In the next theorem, we grant the stochastic comparisons of series systems in the sense of hazard rate order and likelihood ratio order, when the parameters $\alpha_k$ and $\beta_k $ varies.

\begin{theorem}\label{(3.2)}
Consider the setup $X_{k}\thicksim LFR(\alpha_k,\beta_k)$ and $Y_{k}\thicksim LFR(\alpha^*_k,\beta^*_k)$ where $k=1,\ldots,n$. If $\alpha_k\geq\alpha^*_k$ and $\beta_k\geq\beta^*_k$, then $X_{1:n}\leqslant_{hr}Y_{1:n}$.
\end{theorem}
{\bf Proof:} The $HRF$s of $X_{1:n}$ and $Y_{1:n}$ can be written as
\begin{equation}
h_{X_{1:n}}(x)=\sum\limits_{k=1}^n\big(\alpha_k+\beta_k x\big)
\end{equation}
and
\begin{equation}
h_{Y_{1:n}}(x)=\sum\limits_{k=1}^n\big(\alpha^*_k+\beta^*_k x\big),
\end{equation}
respectively. According to Definition \ref{def2.1}(ii), we get $\alpha_k+\beta_k x\geq\alpha^*_k+\beta^*_k x$. This inequality is true for $\alpha_k\geq\alpha^*_k$ and $\beta_k\geq\beta^*_k$, for each $k=1,\ldots,n$, implies that $X_{1:n}\leqslant_{hr}Y_{1:n}$.
\hfill\(\blacksquare\)

We provide Example \ref{example5.2} and Counterexample \ref{counexam5.2} in Section \ref{section5} to support the established result in Theorem \ref{(3.2)}.
{The following corollary is immediate consequence of Theorem \ref{(3.1)}.}
\begin{corollary}
In {Theorem \ref{(3.2)}}, consider $\alpha_k=\alpha^*_k =\alpha$. Then, for fixed $\alpha>0$, we have $X_{1:n}\leqslant_{hr}Y_{1:n}$ under the condition $\beta_k\geq\beta^*_k$.
\end{corollary}

{In the following theorem, we established ordering result between the lifetimes of two series systems in the sense of the likelihood ratio order. In particular, here, we considered heterogeneity occurs in the shape parameters whereas the scale parameters are common and fixed.}
\begin{theorem}\label{(3.3)}
Suppose that the random variables $X_{k}\thicksim LFR(\alpha_k,\beta)$ and $Y_{k}\thicksim LFR(\alpha^*_k,\beta)$. Then, for fixed $\beta>0$, we have
$\alpha_k\geq\alpha^*_k\Rightarrow X_{1:n}\leqslant_{lr}Y_{1:n}$.
\end{theorem}
{\bf Proof:} Let $f_{X_{1:n}}{(x)}$ and $g_{Y_{1:n}}{(x)}$ be the {$PDF$s} of $X_{1:n}$ and $Y_{1:n}$ are given by
\begin{equation}
f_{X_{1:n}}(x)=\exp\left\{-\sum\limits_{k=1}^n\big(\alpha_k x+\frac{\beta}{2}x^2\big)\right\} \sum\limits_{k=1}^n\big(\alpha_k+\beta x\big)
\end{equation}
and
\begin{equation}
g_{Y_{1:n}}(x)=\exp\left\{-\sum\limits_{k=1}^n\big(\alpha^*_k x +\frac{\beta}{2}x^2\big)\right\} \sum\limits_{k=1}^n\big(\alpha^*_k+\beta x\big),
\end{equation}
respectively. Using Part $(iii)$ of Definition \ref{def2.1}, we need to prove $\frac{g_{Y_{1:n}}(x)}{f_{X_{1:n}}(x)}$ is increasing in $x>0$. { Now, the first order partial derivative of $\frac{g_{Y_{1:n}}(x)}{f_{X_{1:n}}(x)}$ with respect to $x$, is obtained as}
\begin{center}
$\frac{n\beta-\Big(\sum\limits_{k=1}^n\big(\alpha^*_k+\beta x\big)\Big)^2}{\sum\limits_{k=1}^n\big(\alpha^*_k+\beta x\big)}-\frac{n\beta-\Big(\sum\limits_{k=1}^n\big(\alpha_k+\beta x\big)\Big)^2}{\sum\limits_{k=1}^n\big(\alpha_k+\beta x\big)}\geq 0,$
\end{center}
under the condition $\alpha_k\geq\alpha^*_k$, {for each $k=1,\ldots,n$ and for fixed $\beta>0$. Thus, the proof is completed.}
\hfill\(\blacksquare\)

In Section \ref{section5}, we provide Example \ref{example5.3} and Counterexample \ref{counexam5.3} to verify the outcome in Theorem \ref{(3.3)}. Next, we enhance certain stochastic comparison results in this section employing the following theorem, which describe the lifetimes of the parallel systems with respect to the usual stochastic order. This ordering result is obtained when the shape and scale parameters are varied for both observations.
\begin{theorem}\label{(3.4)}
For $k=1,\ldots,n$, consider the random variables $X_k\thicksim LFR(\alpha_k,\beta_k)$ and $Y_k\thicksim LFR(\alpha^*_k,\beta^*_k)$. If $\alpha_k\geq\alpha^*_k$ and $\beta_k\geq\beta^*_k$, then $X_{n:n}\leqslant_{st}Y_{n:n}$.
\end{theorem}
{\bf Proof:} The $CDF$s of $X_{n:n}$ and $Y_{n:n}$ are given by
\begin{equation}
F_{X_{n:n}}(x)=\prod\limits_{k=1}^n \Big(1-\exp\left\{-\big(\alpha_k x+\frac {\beta_k}{2}x^2\big)\right\}\Big)
\end{equation}
and
\begin{equation}\label{(3.9)}
G_{Y_{n:n}}(x)=\prod\limits_{k=1}^n \Big(1-\exp\left\{-\big(\alpha^*_k x+\frac {\beta^*_k}{2}x^2\big)\right\}\Big),
\end{equation}
{respectively. By using Part $(i)$ of Definition \ref{def2.1}, we have to show}
\begin{equation}\label{ine:(3.9)}
\prod\limits_{k=1}^n\Big(1-e^{-\big(\alpha_kx+\frac{\beta_k}{2}x^2\big)}\Big)\geq\prod\limits_{k=1}^n\Big(1-e^{-\big(\alpha^*_kx+\frac{\beta^*_k}{2}x^2\big)}\Big).
\end{equation}
{Now, the above inequality given in (\ref{ine:(3.9)}) is satisfied, since the following inequality holds, that is}
\begin{equation}
\Big(1-e^{-\big(\alpha_kx+\frac{\beta_k}{2}x^2\big)}\Big)\geq\Big(1-e^{-\big(\alpha^*_kx+\frac{\beta^*_k}{2}x^2\big)}\Big)
\end{equation}
under the condition $\alpha_k\geq\alpha^*_k$ and $\beta_k\geq\beta^*_k$, which conclude that $X_{n:n}\leqslant_{st}Y_{n:n}$. This completes the proof of the theorem.
\hfill\(\blacksquare\)

Example \ref{example5.4} and Counterexample \ref{counexam5.4} are shown in Section \ref{section5} to validate the established result in Theorem \ref{(3.4)}.

\section{Results based on Transform and Variability Orders}\label{section4}
In this section, we established various ordering results for the lifetime of series systems in the sense of dispersive, star, Lorenz and convex orderings. A random variable {$X\sim LFR(\alpha_k,\beta_k)$, for each $k=1,\ldots,n$}, since its $CDF$ of series system is given in (\ref{(eqn.3.1)}). Now, to obtain the quantile function of $CDF$, set $F(x)=p,~0<p<1$, which is given by
\begin{equation}\label{eq4.1}
1-p=\exp\left\{-\sum\limits_{k=1}^n(\alpha_k x+\frac{\beta_k}{2}x^2)\right\}.
\end{equation}
Now, taking logarithm on both sides of (\ref{eq4.1}), we get
\begin{equation}
\log (1-p)=-\Big(\sum\limits_{k=1}^n(\alpha_k x+\frac{\beta_k}{2}x^2)\Big),
\end{equation}
which is equivalent to
\begin{equation}
\frac{\sum\limits_{k=1}^n\beta_k}{2}x^2+\sum\limits_{k=1}^n\alpha_k x+\log (1-p)=0.
\end{equation}
Now, the above expression is in the form of quadratic equation $ax^2+bx+c=0$, where $a=(\sum_{k=1}^n\beta_k)/2$, $b=\sum_{k=1}^n\alpha_k$, and $c= \log (1-p)$. Therefore, after simplification we obtained the inverse function of $CDF$ $F_{X_{1:n}}(x)$, is given by
\begin{equation}\label{(eqn.(4.2))}
F^{-1}_{X_{1:n}}(p)=\frac{-\sum\limits_{k=1}^n\alpha_k+ \sqrt{(\sum\limits_{k=1}^n\alpha_k)^2-2\sum\limits_{k=1}^n\beta_k \log (1-p)}}{\sum\limits_{k=1}^n\beta_k}.
\end{equation}
Similarly, when $Y\sim LFR(\alpha^*_k,\beta^*_k)$, the quantile function is given by
\begin{equation}\label{(eqn.(4.3))}
F^{-1}_{Y_{1:n}}(p)=\frac{-\sum\limits_{k=1}^n\alpha^*_k+\sqrt{(\sum\limits_{k=1}^n\alpha^*_k)^2-2\sum\limits_{k=1}^n\beta^*_k \log (1-p)}}{\sum\limits_{k=1}^n\beta^*_k}.
\end{equation}
In the following theorem, the dispersive order is established between the lifetimes of two series systems of heterogeneous linear failure rate distributed components, assuming that the distributional shape parameters are connected by inequalities. Here, we considered a common and fixed scale parameter $\beta$ for two series systems.

\begin{theorem}\label{(3.5)}
Consider the setup $X_k\thicksim LFR(\alpha_k,\beta)$ and $Y_k\thicksim LFR(\alpha^*_k,\beta)$ for each $k=1,\ldots,n$. Then, for fixed $\beta>0$, and $\alpha_k\geq\alpha^*_k$, then $X_{1:n}\leqslant_{disp}Y_{1:n}$.
\end{theorem}
{\bf Proof:}
The inverse $CDF$s of $X_{1:n}$ and $Y_{1:n}$ are given in (\ref{(eqn.(4.2))}) and (\ref{(eqn.(4.3))}), respectively, where $0<p<1$. Now, according to Part $(i)$ of Definition \ref{def2.2}, we have to show that $G^{-1}_{Y_{1:n}}(p)-F^{-1}_{X_{1:n}}(p)>0$ is increasing in $p\in\mathbb(0,1)$, {$\frac{\partial}{\partial p}\big(G^{-1}_{Y_{1:n}}(p)-F^{-1}_{X_{1:n}}(p)\big)>0$}. {In doing so,} the {partial} derivative of $G^{-1}_{Y_{1:n}}(p)-F^{-1}_{X_{1:n}}(p)$ with respect to $p$ is {obtained as}
\begin{equation}\label{(3.13)}
\frac{\partial}{\partial p}\Big(G^{-1}_{Y_{1:n}}(p)-F^{-1}_{X_{1:n}}(p)\Big)=\sqrt{\Big(\sum\limits_{k=1}^n \alpha^*_k\Big)^2-2n\beta\log\big(1-p\big)}-\sqrt{\Big(\sum\limits_{k=1}^n \alpha_k\Big)^2-2n\beta\log\big(1-p\big)}.
\end{equation}
After some simplification, we get
\begin{equation}\label{(3.12)}
\Big(\sum\limits_{k=1}^n\alpha^*_k\Big)^2-2n\beta \log\big(1-p\big)<\Big(\sum\limits_{k=1}^n\alpha_k\Big)^2-2n\beta\log\big(1-p\big).
\end{equation}
{Now}, from (\ref{(3.12)}), it is clear that the {partial} derivative of $G^{-1}_{Y_{1:n}}(p)-F^{-1}_{X_{1:n}}(p)$ with respect to $p$ is positive {whenever $\alpha_k\geq\alpha^*_k$}. {Thus,} $G_{Y_{1:n}}^{-1}(p)- F_{X_{1:n}}^{-1}(p)$ is increasing in $p$. {Hence,} under the condition  $\alpha_k\geq\alpha^*_k$,
we can conclude that $X_{1:n}\leq_{disp}Y_{1:n}$ holds.
\hfill\(\blacksquare\)

In Section \ref{section5}, we provide Example \ref{example5.5} and Counterexample \ref{counexam5.5} to verify the established result in Theorem \ref{(3.5)}. In the subsequent result, we investigated the star order between the lifetimes of two series systems $X_{1:n}$ and $Y_{1:n}$ under some sufficient conditions. Here, the shape parameter vectors are varying while the scale parameters are common and fixed for the both systems.

\begin{theorem}\label{(3.6)}
Let $\{X_1,\ldots,X_n\}$ and $\{Y_1,\ldots,Y_n\}$ be two sets of independent random variables where $X_k\thicksim LFR(\alpha_k,\beta)$ and $Y_k\thicksim LFR(\alpha^*_k,\beta)$ for {each} $k=1,\ldots,n$. Then, for fixed $\beta>0$, we have $\alpha^*_k>\alpha_k\Rightarrow X_{1:n}\leqslant_{*}Y_{1:n}$.
\end{theorem}
{\bf Proof:}
Under the given setup, we substitute (\ref{(3.1)}) in (\ref{(eqn.(4.3))}), and obtained
\begin{equation}\label{eqn.(4.5)}
\frac{G^{-1}_{Y_{1:n}}\big({F_{X_{1:n}}(x)}\big)}{x}=\frac{-\sum\limits_{k=1}^n\alpha^*_k+\sqrt{\Big(\sum\limits_{k=1}^n\alpha^*_k\Big)^2+2n\beta x\sum\limits_{k=1}^n\alpha_k +n^2\beta^2x^2}}{n\beta x}.
\end{equation}
{Now}, according to {Part $(ii)$ of} Definition \ref{def2.2}, we need to show that $\frac{G^{-1}_{Y_{1:n}}({F_{X_{1:n}}(x)})}{x}$ is increasing in $x$. {Thus,} to get the desired result, we need to prove  $\frac{\partial}{\partial x}\Big(\frac{G^{-1}_{Y_{1:n}}({F_{X_{1:n}}(x)})}{x}\Big)>0$. {In doing so,} the first {order partial} derivative of ${\frac{G^{-1}_{Y_{1:n}}\big({F_{X_{1:n}}(x)}\big)}{x}}$ with respect to $x$ is obtained as
\begin{align*}
\frac{\partial\big({{G^{-1}_{Y_{1:n}}({F_{X_{1:n}}(x)})}}\big)}{\partial x}&=\sum\limits_{k=1}^n\alpha^*_k \sqrt{\Big(\sum\limits_{k=1}^n\alpha^*_k\Big)^2+2n\beta x\sum\limits_{k=1}^n\alpha_k+\beta^2n^2x^2}-\Big(\Big(\sum\limits_{k=1}^n\alpha^*_k\Big)^2+n\beta x\sum\limits_{k=1}^n\alpha_k\Big)\\
&>0,
\end{align*}
{since $\alpha_k<\alpha^*_k$}. After some simplification, we get, $\alpha_k<\alpha^*_k$ implies that the star order holds between two series systems. This completes the proof of the theorem.
\hfill\(\blacksquare\)

In Section \ref{section5}, we provide Example \ref{example5.6} and Counterexample \ref{counexam5.6} to support the result in Theorem \ref{(3.6)}. The following corollary is a direct consequence of Theorem (\ref{(3.6)}) due to the relationship between the star order and the Lorenz order.
\begin{corollary}
Based on the assumptions and conditions as in Theorem \ref{(3.6)}, we have
\begin{eqnarray*}
\alpha_k<\alpha^*_k\Rightarrow X_{1:n}\leq_{Lorenz}Y_{1:n}.
\end{eqnarray*}	
\end{corollary}
In the next result, the convex transform order between two series systems $X_{1:n}$ and $Y_{1:n}$ has been established. Here, we have considered heterogeneity in the shape parameter vectors while the scale parameter is common and fixed for both series systems.

\begin{theorem}\label{(3.7)}
For $k=1,\ldots,n$, consider the random variables $X_k\thicksim LFR(\alpha_k,\beta)$ and $Y_k\thicksim LFR(\alpha^*_k,\beta)$. Then, for fixed $\beta>0$, we have $\alpha^*_k \geq\alpha_k\Rightarrow X_{1:n}\leqslant_{c}Y_{1:n}$.
\end{theorem}
{\bf Proof:}
Using {Part $(iii)$ of} Definition \ref{def2.2}, we {need} to prove that $G^{-1}_{Y_{1:n}}(F_{X_{1:n}}(x))$ is convex in $x$. For proving its convexity, {it is required to show that
\begin{eqnarray}
\frac{\partial^2G^{-1}_{Y_{1:n}}(F^{-1}_{X_{1:n}}(x))}{\partial x^2}\geq 0.
\end{eqnarray}
Now, for doing so,} differentiating {partially} twice (\ref{eqn.(4.5)}) with respect to $x$, we obtained
\begin{align*}
\frac{\partial^2\big({{G^{-1}_{Y_{1:n}}({F_{X_{1:n}}(x)})}}\big)}{\partial x^2}&=\Big(\sum\limits_{k=1}^n\alpha^*_k\Big)^2+2n\beta\Big(\sum\limits_{k=1}^n\alpha_kx+\frac{\beta}{2}x^2\Big)-\Big(\sum\limits_{k=1}^n\alpha_k+\beta x\Big)^2\\
&\geq 0,
\end{align*}
under the condition $\alpha^*_k\geq\alpha_k$, {which} conclude that $X_{1:n}\leqslant_{c}Y_{1:n}$. {This completes the proof of the theorem.}
\hfill\(\blacksquare\)

Section \ref{section5} includes Example \ref{example5.7} and Counterexample \ref{counexam5.7} to validate the established result in Theorem \ref{(3.7)}.

\section{Numerical Examples and Counterexamples}\label{section5}
\begin{figure}[H]
\begin{center}
\subfigure[]{\includegraphics[width=3.2in,height=2.5in]{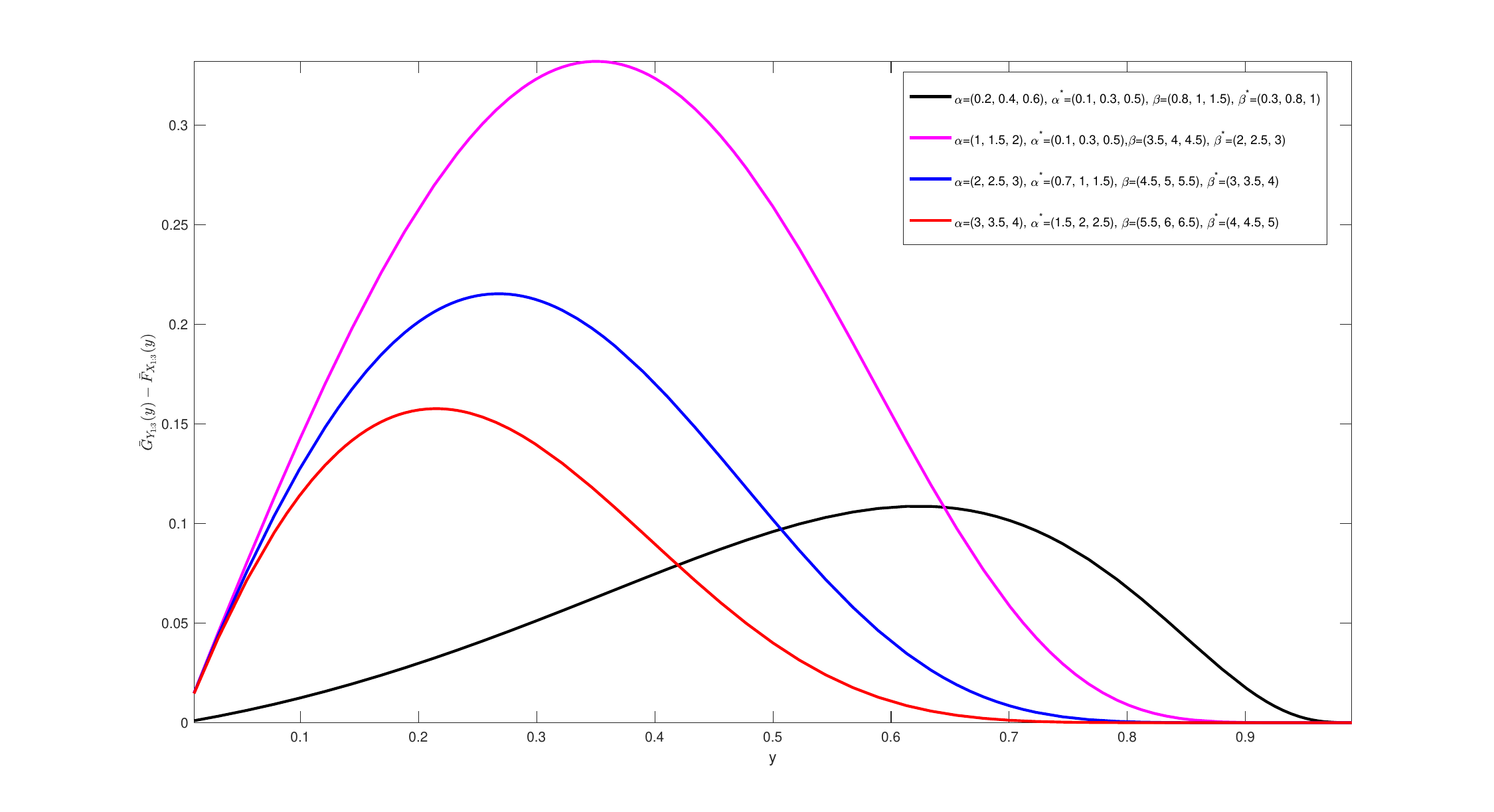}}
\subfigure[]{\includegraphics[width=3.2in,height=2.5in]{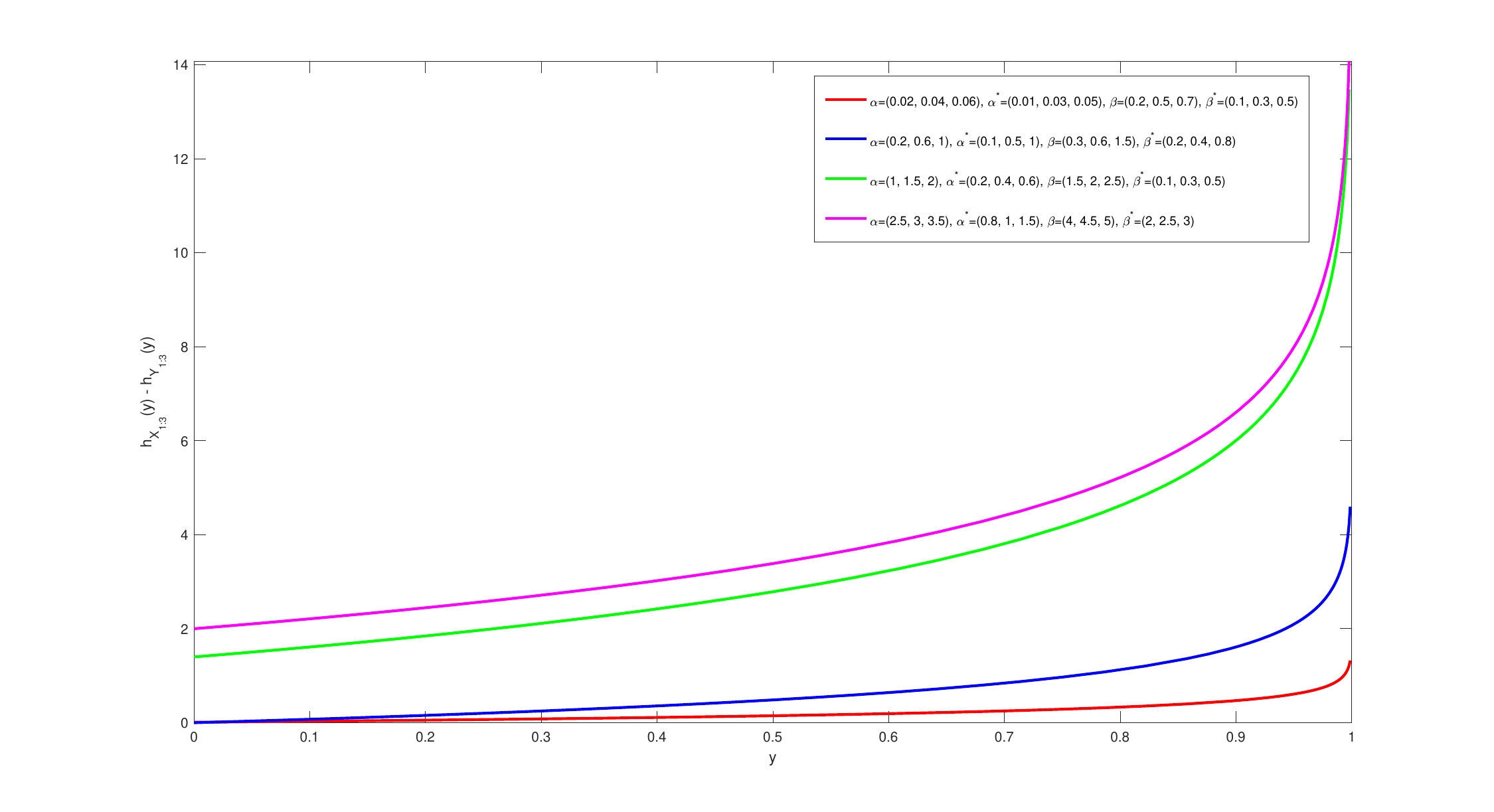}}
\caption{{(a) Plots of the difference $\bar{G}_{Y_{1:3}}(y)-\bar{F}_{X_{1:3}}(y)$ as in Example \ref{example5.1}}. (b) Plots {of the difference $h_{X_{1:3}}(y)-h_{Y_{1:3}}(y)$ as in Example \ref{example5.2}.}}
\label{(1)}
\end{center}
\end{figure}

Before exploring examples, let us discuss the domain of the function where $x\in\mathbb{R}^+$ (the set of all positive real numbers), and it is challenging to graph it over this domain, one possible approach is to map the variable $x$ to new variable $y$ defined in the range $(0.01,0.99)$. This can be achieved using a transformation function, such as: $x=\log(\frac{1}{1-y})$, for $y\in(0,1)$. Here, this transformation allows us to analyze and graph the function more effectively. By restricting the domain in this way, we can focus on a manageable range of values and gain insights into the behavior of the function. This approach helps to compress an infinite domain $(\mathbb{R}^+)$ into a finite interval $(0,1)$ for easier graphing while preserving the function's behavior. In this section, we demonstrated various numerical examples {and counterexamples} for validation of the {established} ordering results obtained in Sections  \ref{section3} and \ref{section4}, respectively. The following is an example of the verification of Theorem \ref{(3.1)} and also present a counterexample that highlight the importance of the conditions outlined in the Theorem \ref{(3.1)} for accurate result.
\begin{exmp}\label{example5.1}
Let $\{X_1,X_2,X_3\}$ and $\{Y_1,Y_2,Y_3\}$ be the two sets of independent random variables with $X_k\thicksim\mathcal{LFR}(\alpha_k,\beta_k)$ and $Y_k\thicksim\mathcal{LFR}(\alpha^*_k,\beta^*_k)$ for $k=1,2,3$, respectively. It can be easily verified that the conditions of Theorem \ref{(3.1)} hold. Plot the whole of distribution function curves of $X_{1:3}$ and $Y_{1:3}$ are display in Figure \ref{(1)}(a), which confirmed that $X_{1:3}\leqslant_{st}Y_{1:3}$.
\end{exmp}
{The next counterexample shows the importance of the sufficient conditions ``$\alpha_k\leq\alpha^*_k$'' and ``$\beta_k\leq\beta^*_k$'' to establish the usual stochastic ordering between two series systems in Theorem \ref{(3.1)}.}

\begin{figure}[H]
\begin{center}
\subfigure[]{\includegraphics[width=2in,height=2in]{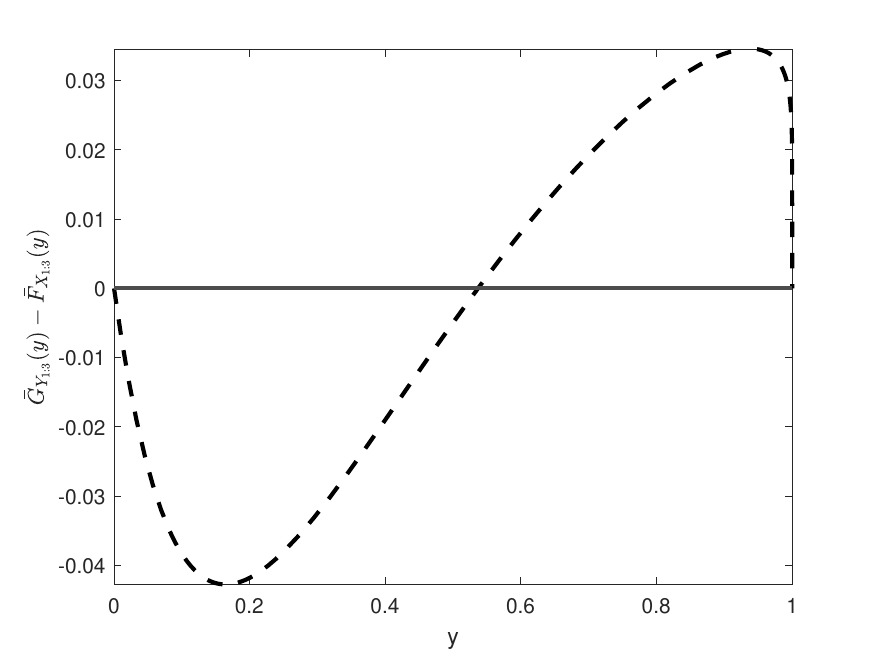}}
\subfigure[]{\includegraphics[width=2in,height=2in]{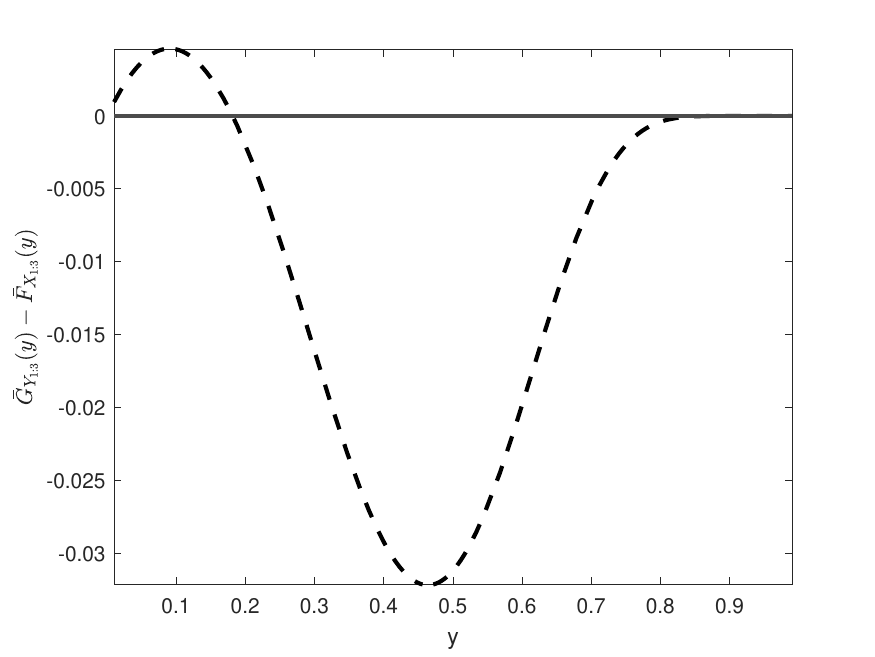}}
\subfigure[]{\includegraphics[width=2in,height=2in]{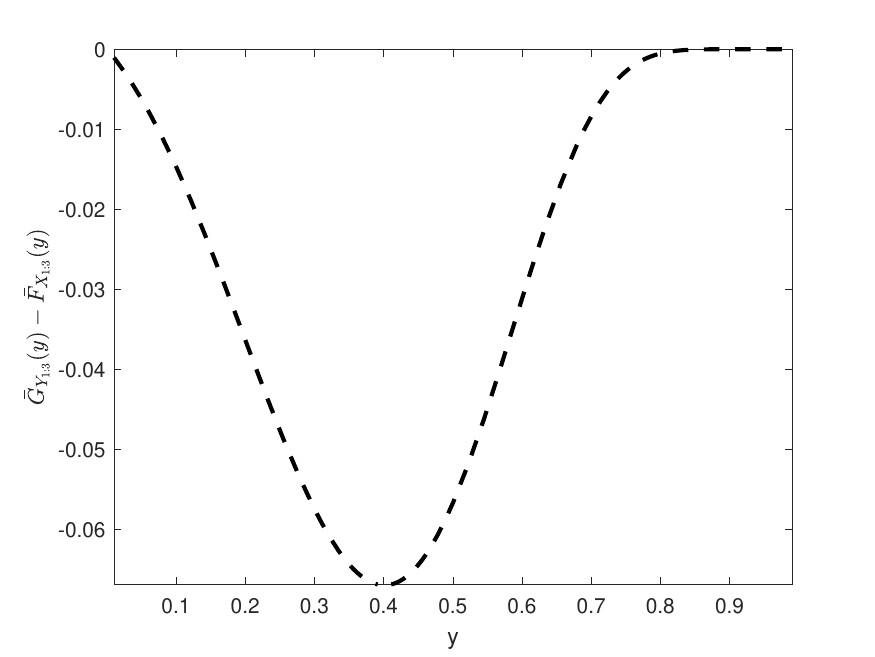}}
\caption{Difference of {the} survival functions of ${Y_{1:3}}$ and ${X_{1:3}}$. $(a)$: represents the difference {between two} survival functions is not always non-negative in $y$ under the conditions $\alpha_k<\alpha^*_k$ {and} $\beta_k>\beta^*_k$. $(b)$: represents the difference {between two} survival functions under the conditions $\alpha_k>\alpha^*_k$ {and} $\beta_k<\beta^*_k$ which is non-negative. $(c)$: represents the curve of $\bar{G}_{Y_{1:3}}(y)-\bar{F}_{X_{1:3}}(y)$ under the conditions $\alpha_k<\alpha^*_k$ {and} $\beta_k<\beta^*_k$ is negative.}
\label{(2)}
\end{center}
\end{figure}

\begin{counterexmp}\label{counexam5.1}
Assume $\alpha_k\leq\alpha^*_k$ {and} $\beta_k\leq\beta^*_k$. Clearly, {the} assumptions made in Theorem \ref{(3.1)} are not satisfied. As shown in Figure \ref{(2)}, the difference between the survival functions $\bar{G}_{Y_{1:3}}(y)$ and $\bar{F}_{X_{1:3}}(y)$ is not always non-negative in $y$. {Thus, ${X_{1:3}}\nleq_{st}{Y_{1:3}}$.}
\end{counterexmp}

{The following example provides an illustration of the result in Theorem \ref{(3.2)}.}
\begin{exmp}\label{example5.2}
Assume that $\alpha_k=(0.02,0.04,0.06)$, $\alpha^*_k=(0.01,0.03,0.05)$, $\beta_k=(0.2,0.5,0.7)$, and $\beta^*_k=(0.1,0.3,0.5)$ {for each} $k=1,2,3$. Clearly, the assumptions made in Theorem \ref{(3.2)} are satisfied for all $k>0$. Taking these numerical values of the parameters, {the difference between} $h_{X_{1:3}}(y)$ and $h_{Y_{1:3}}(y)$ is plotted in Figure \ref{(1)}(b), validating the result in Theorem \ref{(3.2)}.
\end{exmp}

{The following counterexample shows that the desired ordering result in Theorem \ref{(3.2)} does not hold if $\alpha_k\ngeq\alpha^*_k$ and $\beta_k\ngeq\beta^*_k$.}

\begin{figure}[H]
\begin{center}
\subfigure[]{\includegraphics[width=2in,height=2in]{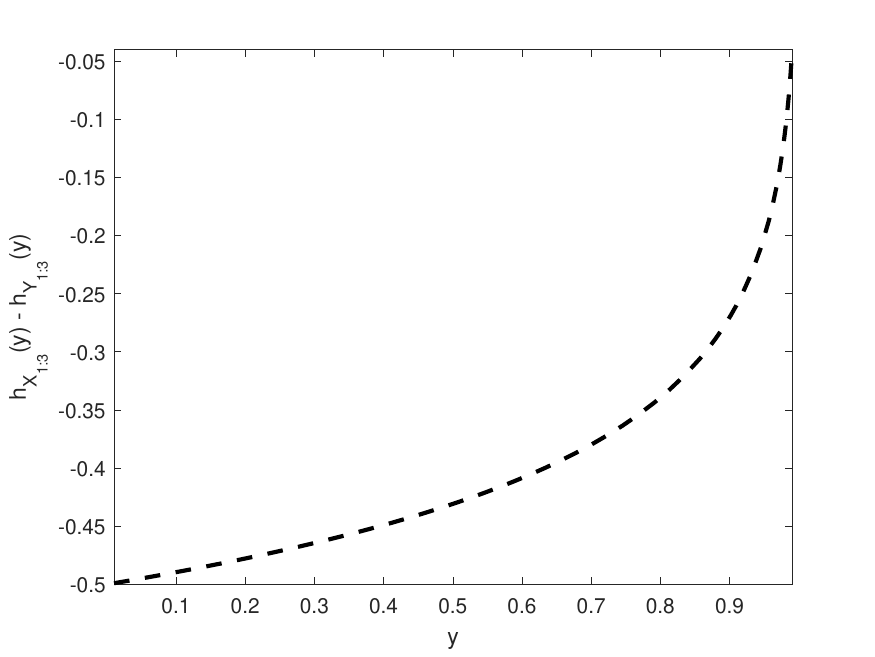}}
\subfigure[]{\includegraphics[width=2in,height=2in]{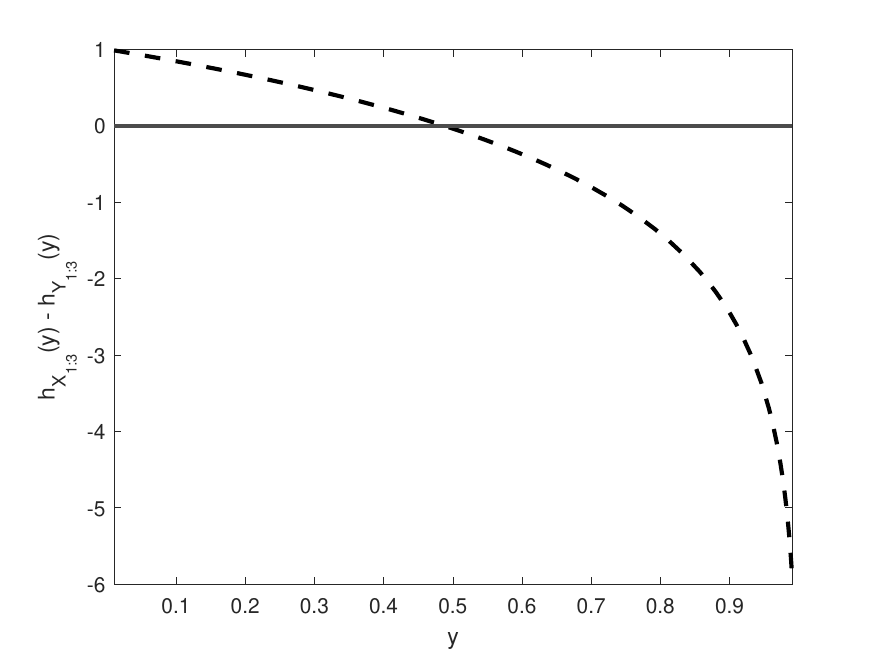}}
\subfigure[]{\includegraphics[width=2in,height=2in]{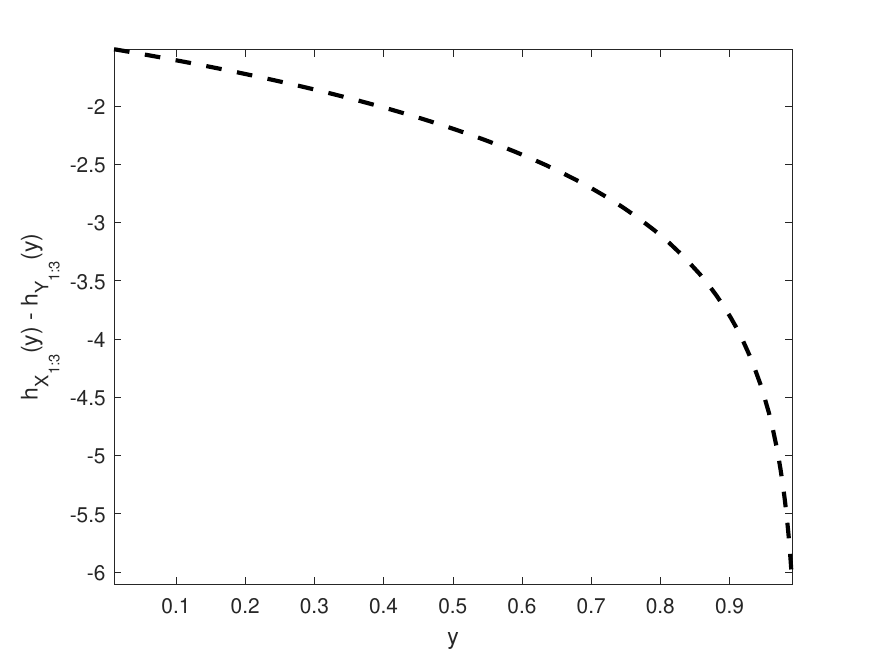}}
\caption{Plots represent the difference {between two} hazard rate functions of ${X_{1:3}}$ and ${Y_{1:3}}$. $(a)$: represents the difference {between two} hazard rate functions which is not always non-negative in $y$, when the conditions $\alpha_k<\alpha^*_k$ {and} $\beta_k>\beta^*_k$ holds. $(b)$: represents the difference {between two} hazard rate functions which is not always non-negative in $y$, {while} the conditions $\alpha_k>\alpha^*_k$ {and} $\beta_k<\beta^*_k$ holds. $(c)$: represents the difference {between two} hazard rate functions which is negative in $y$, when the conditions $\alpha_k<\alpha^*_k$ {and} $\beta_k<\beta^*_k$ holds.}
\label{(4)}
\end{center}
\end{figure}

\begin{counterexmp}\label{counexam5.2}
Let $\{X_1,X_2,X_3\}$ and $\{Y_1,Y_2,Y_3\}$ be the two sets of independent random observations such that $X_k\sim\mathcal{LFR}(\alpha_k,\beta_k)$ and $Y_k\sim\mathcal{LFR}(\alpha^*_k,\beta^*_k)$, {for each $k=1,2,3$.} As mentioned in Figure \ref{(4)}, the difference between {the} hazard rate functions $h_{X_{1:3}}(y)$ and $h_{Y_{1:3}}(y)$ is not always non-negative in $y$, except the conditions $\alpha_k\ngeq\alpha^*_k$ {and} $\beta_k\ngeq\beta^*_k$. {Here,} all other conditions of Theorem \ref{(3.2)} are satisfied. {Clearly, the difference take negative as well as positive values, which} means that $X_{1:3}\nleq_{hr}Y_{1:3}$. So, the result {established} in Theorem \ref{(3.2)} {does not} hold.
\end{counterexmp}
{In order to justify Theorem \ref{(3.3)}, an example is provided.}

\begin{figure}[H]
\begin{center}
\subfigure[]{\includegraphics[width=3.2in,height=2.5in]{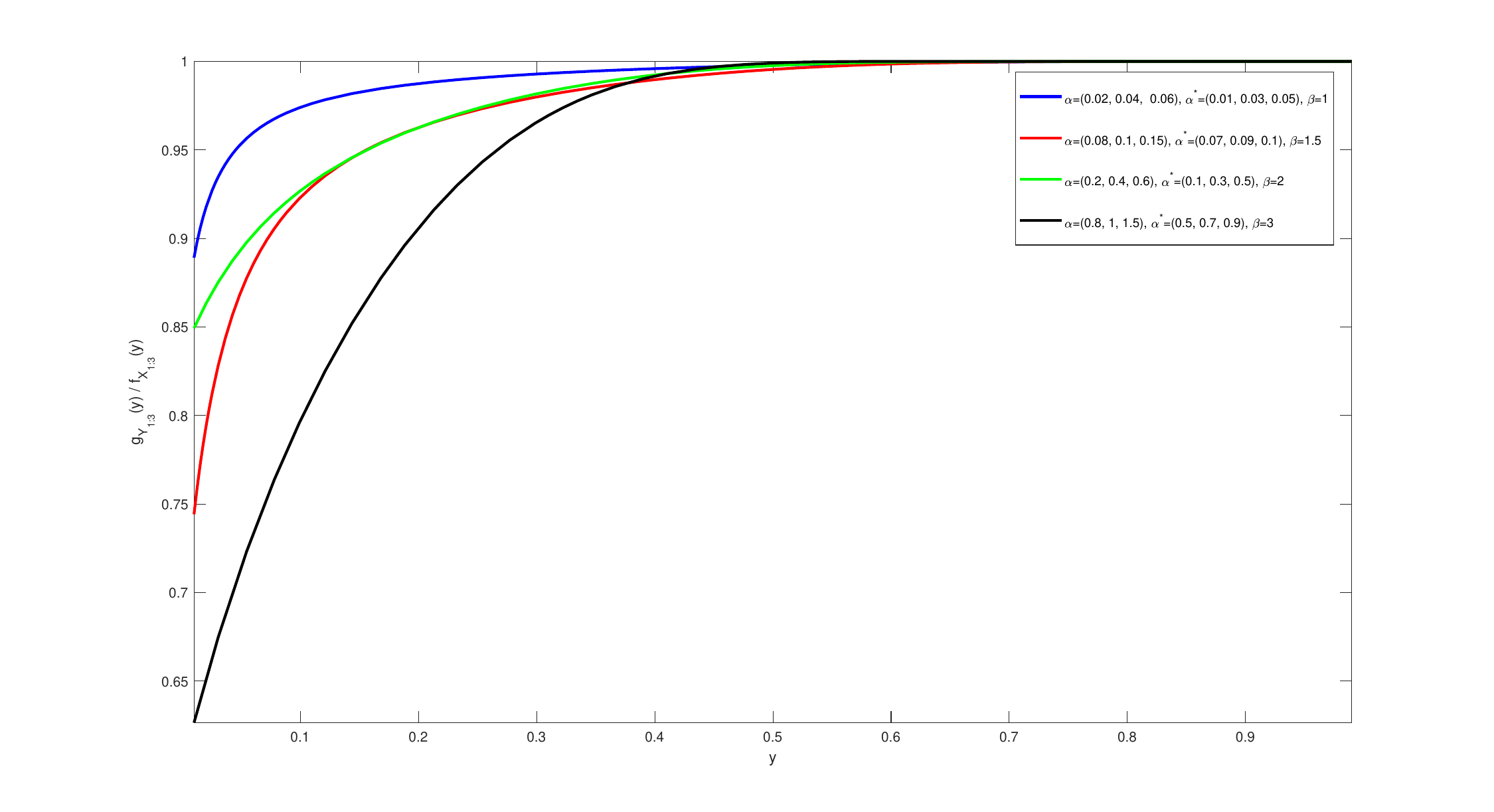}}
\subfigure[]{\includegraphics[width=3.2in,height=2.5in]{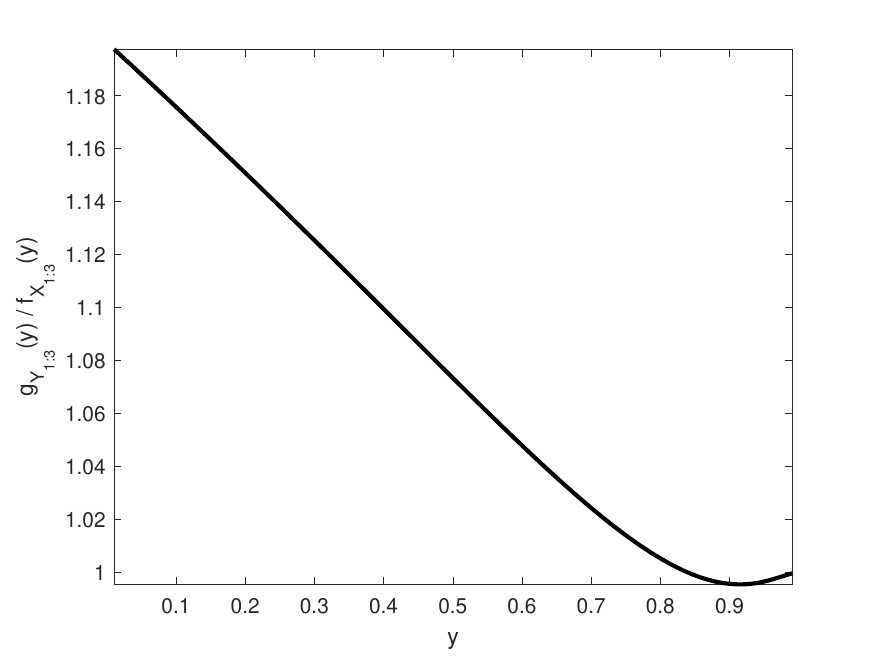}}
\caption{ $(a)$ Plots of the ratio of the {$PDF$s $g_{Y_{1:3}}(y)$ and $f_{X_{1:3}}(y)$ as in Example \ref{example5.3}.} $(b)$ Plot represents the ratio of the {$PDF$s} $g_{Y_{1:3}}(y)$ and $f_{X_{1:3}}(y)$ as in Counterexample \ref{counexam5.3}.}
\label{(5)}
\end{center}
\end{figure}

\begin{exmp}\label{example5.3}
Let $\{X_1,X_2,X_3\}$ and $\{Y_1,Y_2,Y_3\}$ be two sets of independent random variables with $X_k\thicksim\mathcal{LFR}(\alpha_k,\beta)$ and $Y_k\thicksim\mathcal{LFR}(\alpha^*_k,\beta)$ for $k=1,2,3$, respectively. Assume $(\alpha_1,\alpha_2,\alpha_3)=(0.02,0.04,0.06)$, $(\alpha^*_1,\alpha^*_2,\alpha^*_3)=(0.01,0.03,0.05)$, $n=3$, and $\beta=1$. Clearly, {all the} assumptions made in Theorem \ref{(3.3)} are satisfied. {Now,} using the numerical values {of the parameters}, we plot the graphs of $\frac{g_{Y_{1:3}}(y)}{f_{X_{1:3}}(y)}$ {in Figure \ref{(5)}}(a), {which} is increasing in $y$. {From,} Figure \ref{(5)}(a), {it is clear that $X_{1:3}\leq_{lr}Y_{1:3}$}, validating the result of likelihood ratio order in Theorem \ref{(3.3)}.
\end{exmp}

{We now present a counterexample to emphasize that the condition ``$\alpha_k\geq\alpha^*_k$'' is required for the result in Theorem \ref{(3.3)}.}

\begin{counterexmp}\label{counexam5.3}
Assume that $(\alpha_1,\alpha_2,\alpha_3)=(0.1,0.3,0.5)$, $(\alpha^*_1,\alpha^*_2,\alpha^*_3)=(0.2,0.4,0.6)$, and $\beta=0.1$. Clearly, except {the condition} $\alpha_k\ngeq\alpha^*_k$, all {the} other conditions of Theorem \ref{(3.3)} are satisfied. Based on the numerical values of the parameters, the ratio of the {$PDF$s} $g_{Y_{1:3}}(y)$ and $f_{X_{1:3}}(y)$ is given in Figure \ref{(5)}(b). {From this figure,} we can conclude that {the} ratio is a non-monotone function in $y>0$. Thus, the likelihood ratio order between the {series systems} $X_{1:3}$ and $Y_{1:3}$ in Theorem \ref{(3.3)} does not hold.
\end{counterexmp}

{The following example illustrates Theorem \ref{(3.4)}, for $n=3$.}

\begin{figure}[H]
\begin{center}
\subfigure[]{\includegraphics[width=3.2in,height=2.5in]{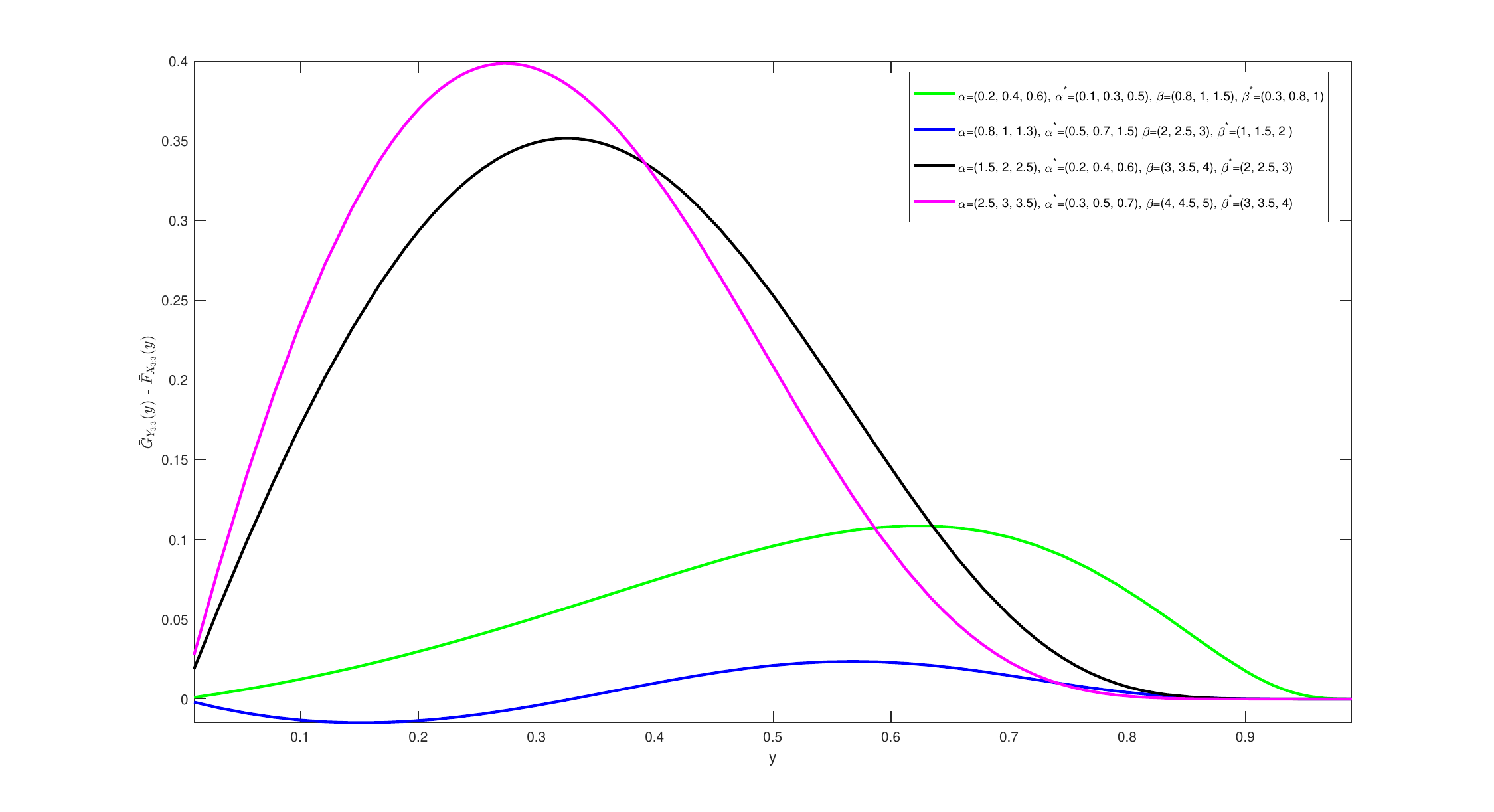}}
\subfigure[]{\includegraphics[width=3.2in,height=2.5in]{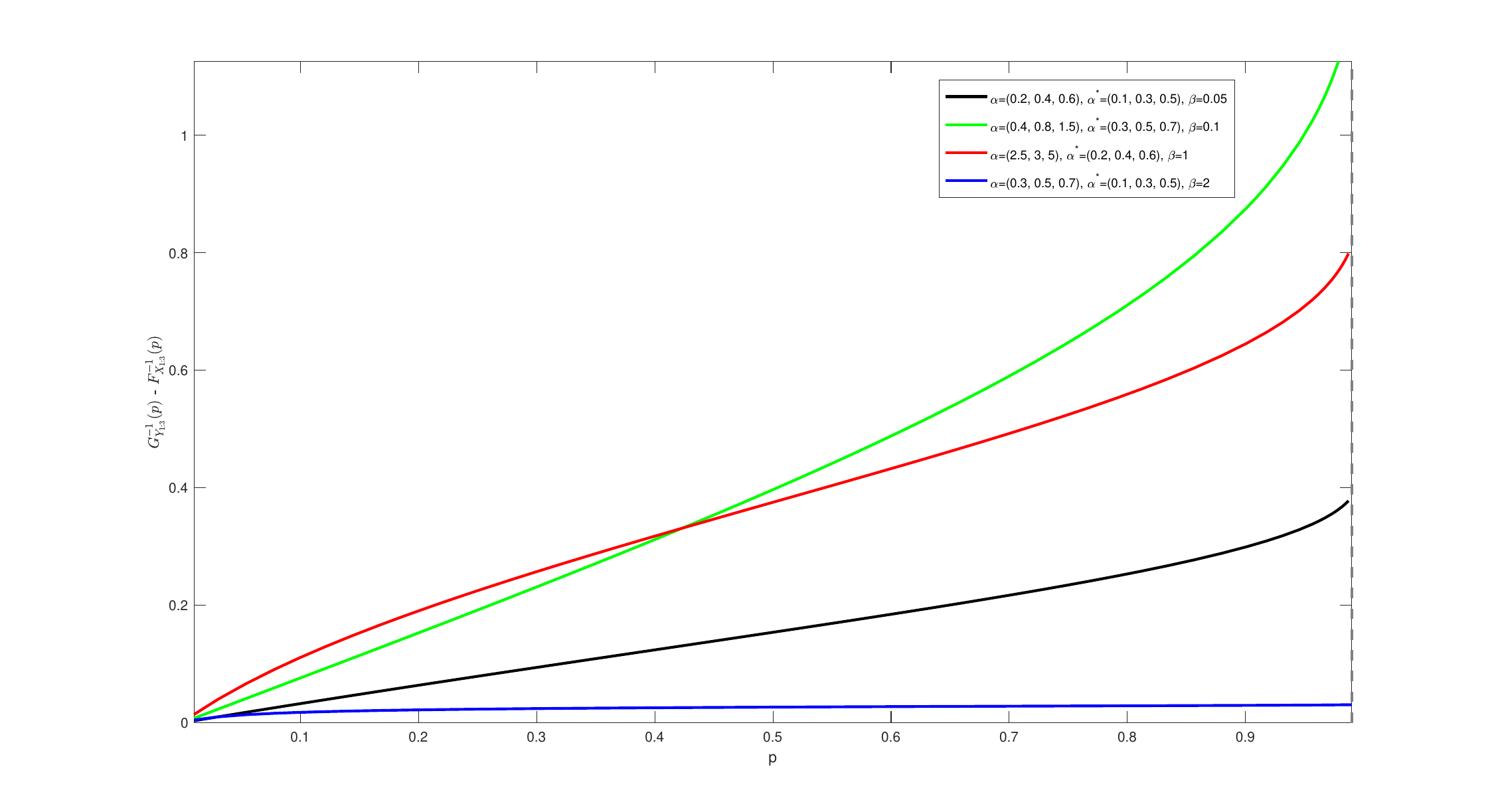}}
\caption{{$(a)$ Plots of the difference $\bar{G}_{Y_{3:3}}(y)-\bar{F}_{X_{3:3}}(y)$ as in Example \ref{example5.4}. $(b)$ Plots of the difference $G^{-1}_{Y_{1:3}}(p)-F^{-1}_{X_{1:3}}(p)$ as in Example \ref{example5.5}.}}
\label{(7)}
\end{center}
\end{figure}

\begin{exmp}\label{example5.4}
Let $\{X_1,X_2,X_3\}$ and $\{Y_1,Y_2,Y_3\}$ be two sets of independent random variables with $X_k\thicksim\mathcal{LFR}(\alpha_k,\beta_k)$ and $Y_k\thicksim\mathcal{LFR}(\alpha^*_k,\beta^*_k)$ for $k=1,2,3$, respectively. In addition, suppose that $(\alpha_1,\alpha_2,\alpha_3)=(0.2,0.4,0.6)$, $(\alpha^*_1,\alpha^*_2,\alpha^*_3)=(0.1,0.3,0.5)$, $(\beta_1,\beta_2,\beta_3)=(0.8,1,1.5)$, and $(\beta^*_1,\beta^*_2,\beta^*_3)=(0.3,0.8,1)$. It can be easily verified that the conditions of Theorem \ref{(3.4)} hold. {Now, we plot the difference between the $SF$s of} $Y_{3:3}$ and $X_{3:3}$ in Figure \ref{(7)}(a), which confirmed that $X_{3:3}\leqslant_{st}Y_{3:3}$.
\end{exmp}
{Next, present a counterexample to show that the sufficient conditions ``$\alpha_k\geq\alpha^*_k$'' and ``$\beta_k\geq\beta^*_k$'' are necessary for establishing the usual stochastic order in Theorem \ref{(3.4)}.}

\begin{figure}[H]
\begin{center}
\subfigure[]{\includegraphics[width=2in,height=2in]{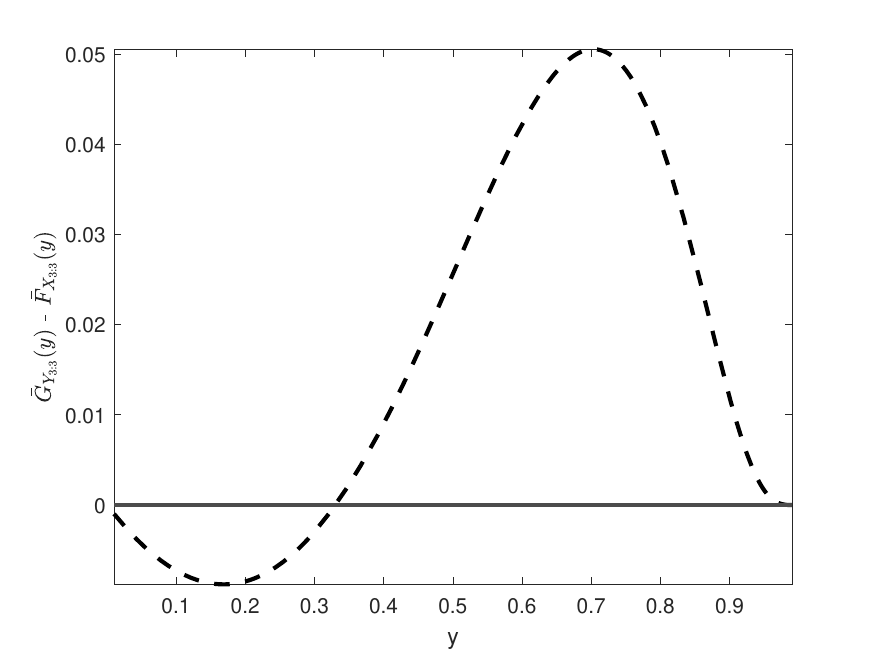}}
\subfigure[]{\includegraphics[width=2in,height=2in]{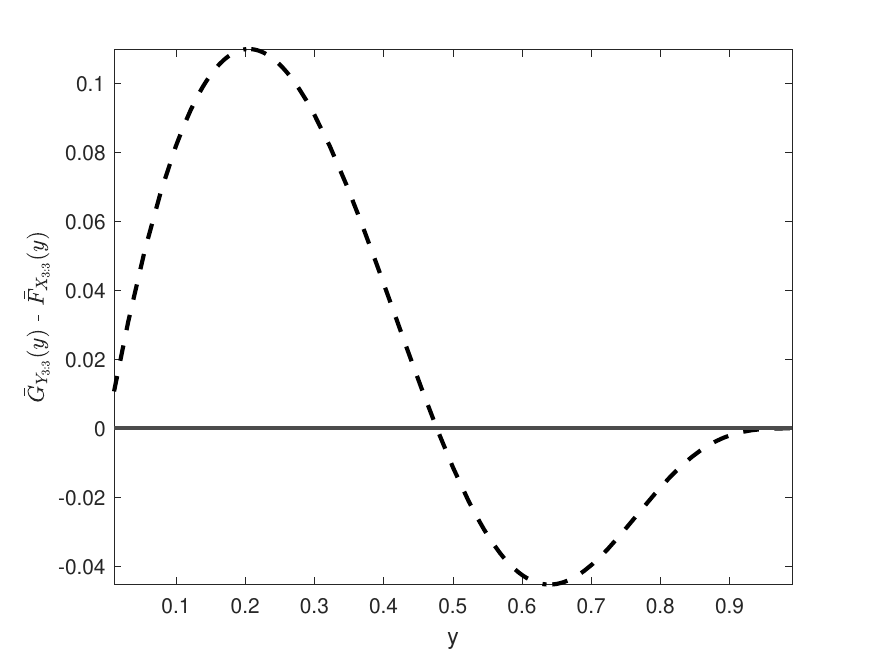}}
\subfigure[]{\includegraphics[width=2in,height=2in]{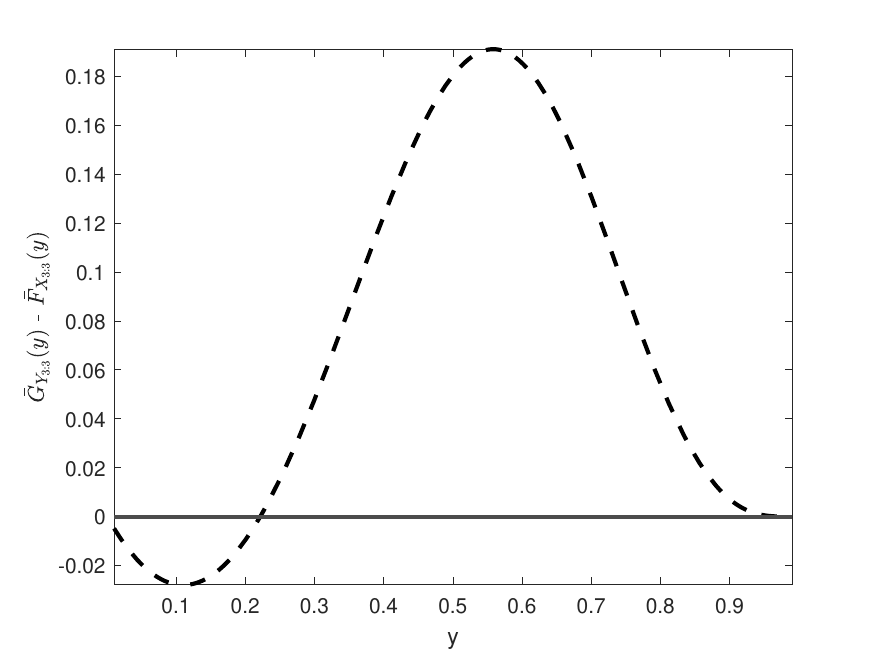}}
\caption{Plots represent the difference {between two $SF$s} of $Y_{3:3}$ and $X_{3:3}$ which is not always non-negative $(a)$: exist, when the conditions $\alpha_k<\alpha^*_k$ {and} $\beta_k<\beta^*_k$ holds. $(b)$: exist, when the conditions $\alpha_k>\alpha^*_k$ {and} $\beta_k<\beta^*_k$ {holds}. $(c)$: exist, when the conditions $\alpha_k<\alpha^*_k$ {and} $\beta_k<\beta^*_k$ {holds}.}
\label{(8)}
\end{center}
\end{figure}

\begin{counterexmp}\label{counexam5.4}
Let $X_k\thicksim\mathcal{LFR}(\alpha_k,\beta_k)$ and $Y_k\thicksim\mathcal{LFR}(\alpha^*_k,\beta^*_k)$ for $k=1,2,3$, respectively. The curve of {the difference $\bar{G}_{Y_{3:3}}(y)-\bar{F}_{X_{3:3}}(y)$} is displayed in Figure \ref{(8)}, which is not always non-negative in $y$ and this means that $X_{3:3}\nleq_{st}Y_{3:3}$. {Thus,} the result {established} in Theorem \ref{(3.4)} {does not} hold.
\end{counterexmp}

{Now, we provide some examples and counterexamples to illustrate the transform and variability orders comparing between two series systems. In order to justify Theorem \ref{(3.5)}, an example is provided.}
\begin{exmp}\label{example5.5}
Let $\{X_1,X_2,X_3\}$ and $\{Y_1,Y_2,Y_3\}$ be two sets of independent random observations such that $X_k\thicksim\mathcal{LFR}(\alpha_k,\beta)$ and $Y_k\thicksim \mathcal{LFR}(\alpha^*_k,\beta)$, $k=1,2,3$. {Under this setting, we have plotted the graph of $G^{-1}_{Y_{1:3}}(p)-F^{-1}_{X_{1:3}}(p)$ in Figure \ref{(7)}(b)}, for $\beta=0.05$, $n=3$, $(\alpha_1,\alpha_2,\alpha_3)=(0.2,0.4,0.6)$, and $(\alpha^*_1,\alpha^*_2,\alpha^*_3)=(0.1,0.3,0.5)$. {From Figure \ref{(7)}(b), it is obvious that $X_{1:3}\leq_{disp}Y_{1:3}$ holds, illustrating the result stated in Theorem \ref{(3.5)}.}
\end{exmp}

\begin{figure}[H]
\begin{center}
\subfigure[]{\includegraphics[width=3.2in,height=2.5in]{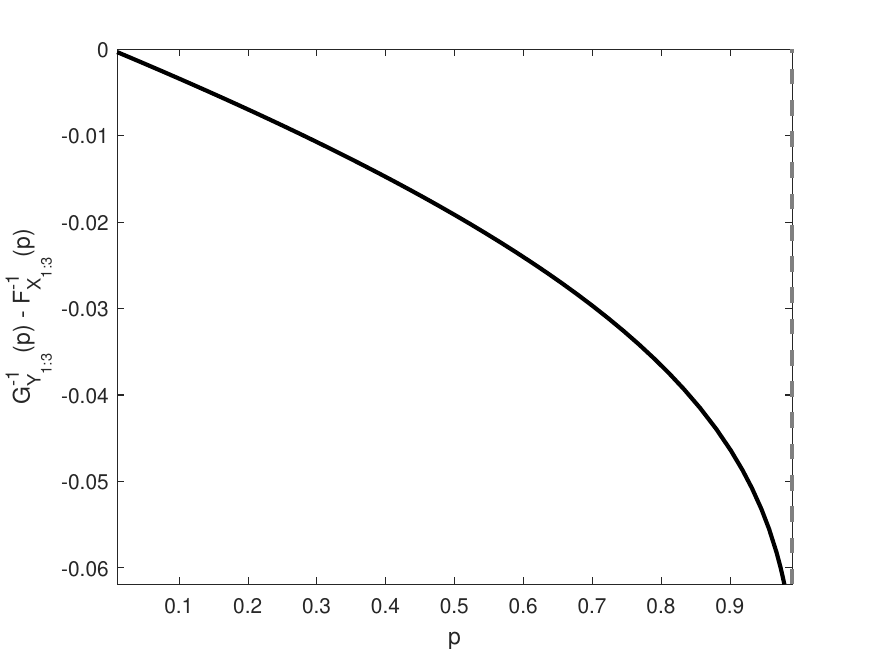}}
\subfigure[]{\includegraphics[width=3.2in,height=2.5in]{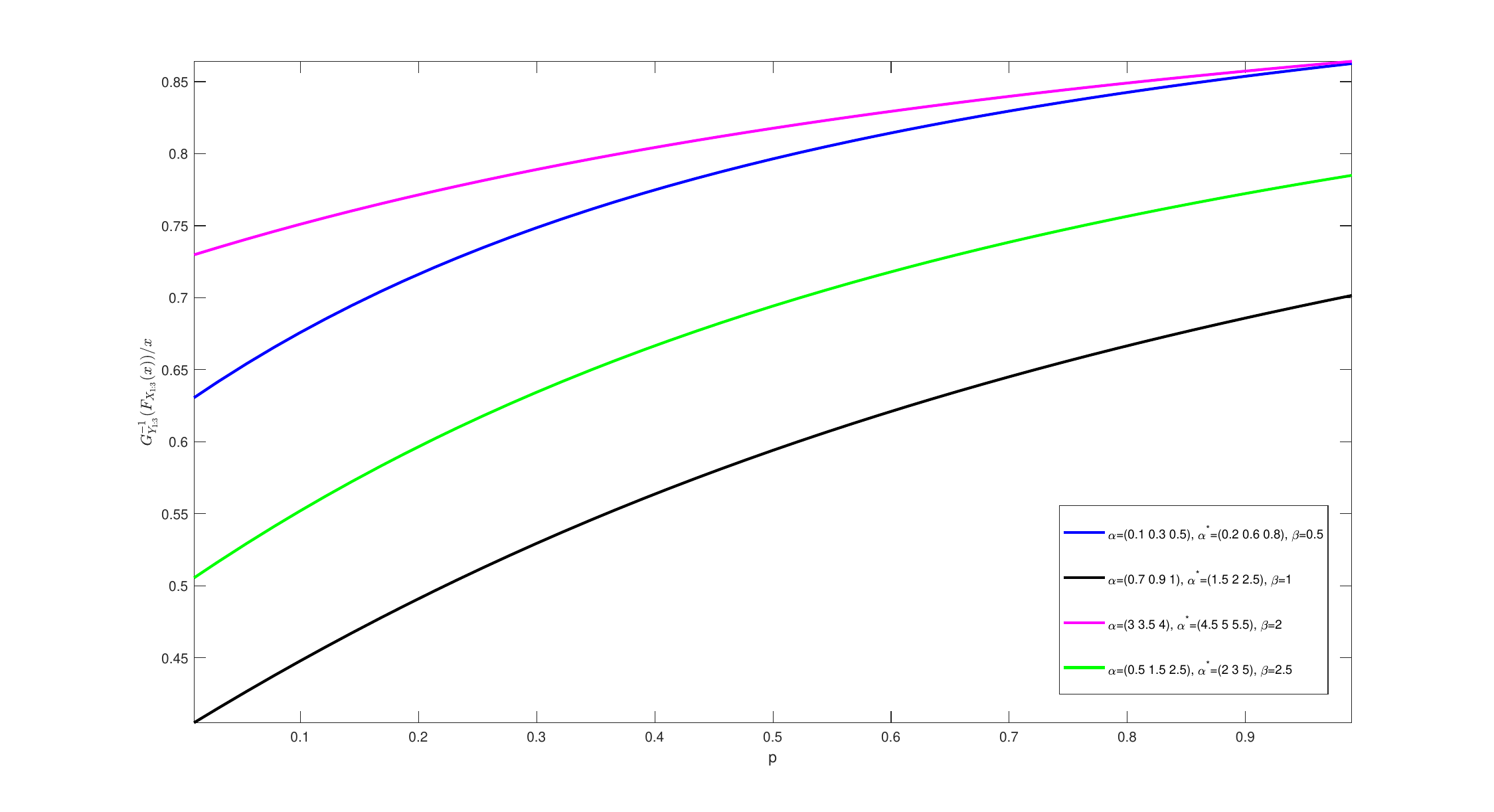}}
\caption{{$(a)$ Plot of the difference $G^{-1}_{Y_{1:3}}(p)-F^{-1}_{X_{1:3}}(p)$ as in Counterexample \ref{counexam5.5}. $(b)$ Plots of $G^{-1}_{Y_{1:3}}\big(F_{X_{1:3}}(x)\big)/x$ as in Example \ref{example5.6}.}}
\label{(10)}
\end{center}
\end{figure}

{The following counterexample illustrates that the result in Theorem \ref{(3.5)} does not hold if ``$\alpha_k\ngeq\alpha^*_k$''.}
\begin{counterexmp}\label{counexam5.5}
Assume that $n=3$, $\beta=2$, $(\alpha_1,\alpha_2,\alpha_3)=(1,3,5)$, and $(\alpha^*_1,\alpha^*_2,\alpha^*_3)=(2,4,6)$. Here, {it is obvious that} $(\alpha_1,\alpha_2,\alpha_3)\ngeq(\alpha^*_1,\alpha^*_2,\alpha^*_3)$, while other conditions are satisfied {in Theorem \ref{(3.5)}}. Now, the curve of {the difference} $G^{-1}_{Y_{1:3}}(p)-F^{-1}_{X_{1:3}}(p)$ is plotted in Figure \ref{(10)}(a). {From figure, we see that} the curve take negative {values}, which means that the desired dispersive order in Theorem \ref{(3.5)} does not hold.
\end{counterexmp}

{Theorem \ref{(3.6)} is illustrated in the next example.}
\begin{exmp}\label{example5.6}
Consider $n=3$, $\beta=0.5$, $(\alpha_1,\alpha_2,\alpha_3)=(0.1,0.3,0.5)$, and $(\alpha^*_1,\alpha^*_2,\alpha^*_3)=(0.2,0.6,0.8)$. Clearly, {all} the assumptions made in Theorem \ref{(3.6)} are satisfied. {Now,} using the numerical values {of the parameters}, we plot the graphs of {$G^{-1}_{Y_{1:3}}\big(F_{X_{1:3}}(x)\big)/x$} in Figure \ref{(10)}(b), validating the {star} order in Theorem \ref{(3.6)}.
\end{exmp}

\begin{figure}[H]
\begin{center}
\subfigure{\includegraphics[width=3.2in,height=2.5in]{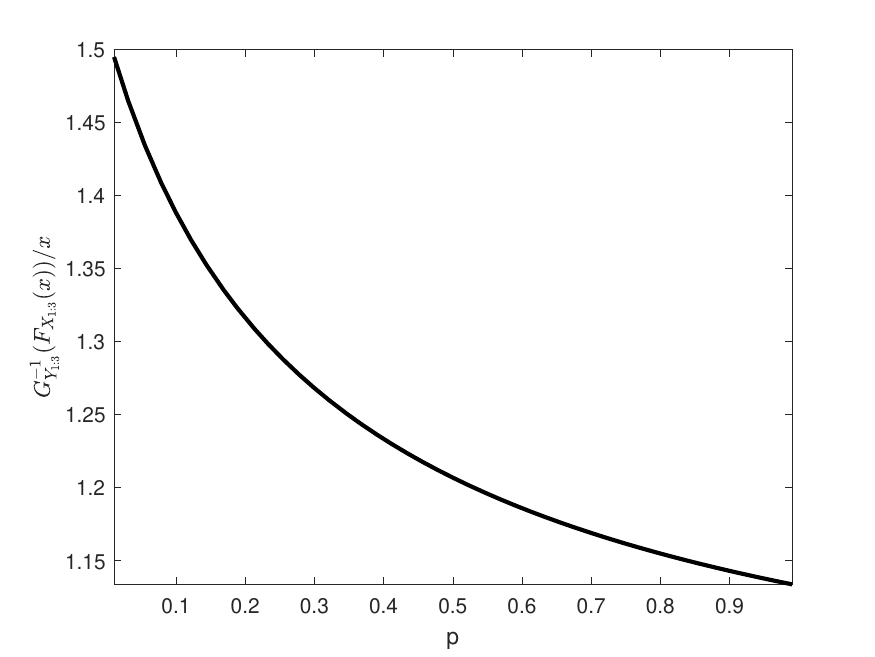}}
\subfigure{\includegraphics[width=3.2in,height=2.5in]{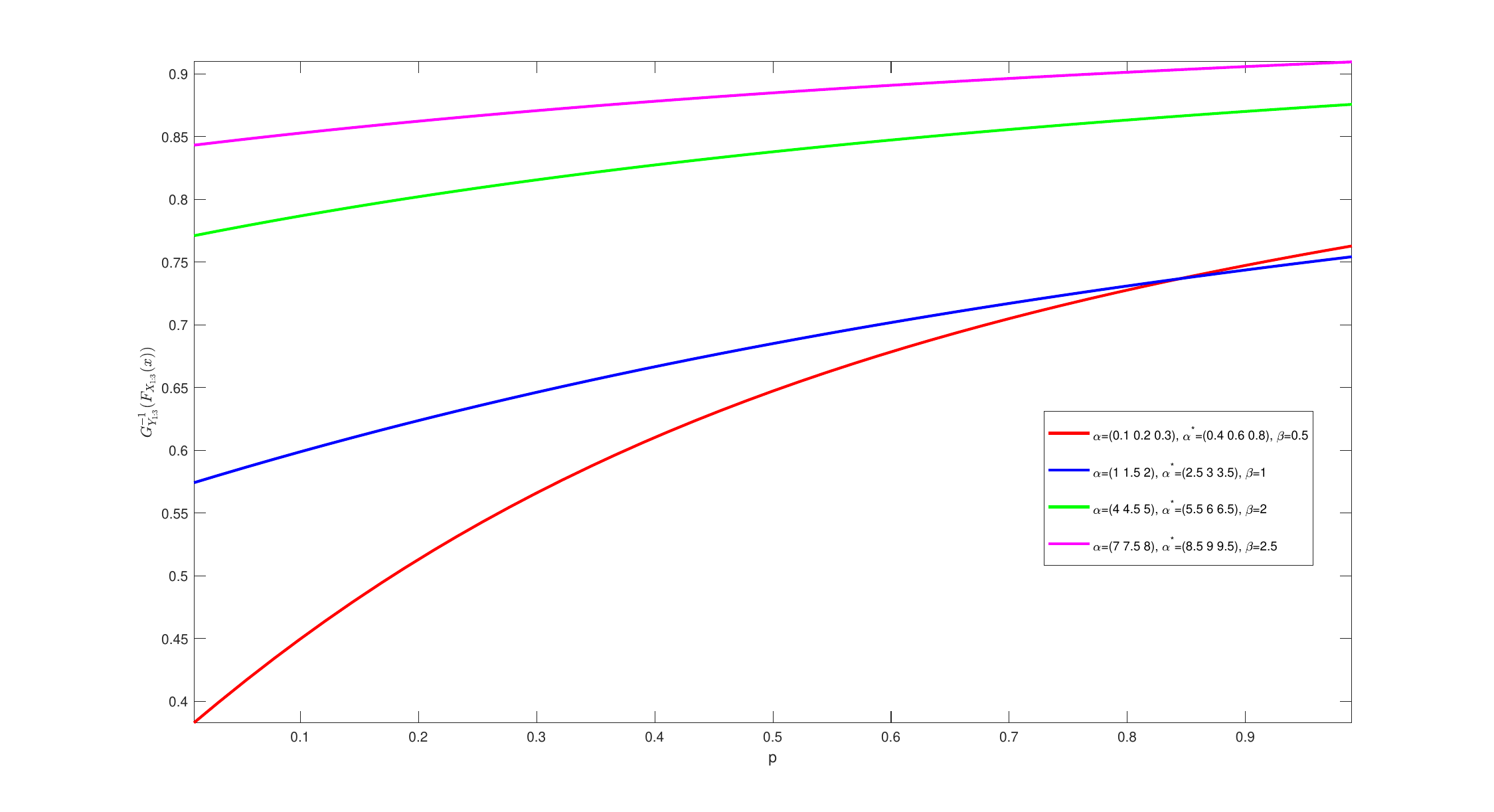}}
\caption{{$(a)$ Plot of $G^{-1}_{Y_{1:3}}\big(F_{X_{1:3}}(x)\big)/x$ as in Counterexample \ref{counexam5.6}. $(b)$ Plots of $G^{-1}_{Y_{1:3}}\big(F_{X_{1:3}}(x)\big)$ as in Example \ref{example5.7}.}}
\label{(12)}
\end{center}
\end{figure}

{In the following counterexample, we illustrate that the condition ``$\alpha^*_k>\alpha_k$'' plays a significant role in establishing the star order between two lifetimes of the series systems.}
\begin{counterexmp}\label{counexam5.6}
Let $n=3$, $\beta=4$, $(\alpha_1,\alpha_2,\alpha_3)=(6.2,6.6,6.8)$, and $(\alpha^*_1,\alpha^*_2,\alpha^*_3)=(4.1,4.3,4.5)$. Clearly, $(6.2,6.6,6.8)\nleq(4.1,4.3,4.5)$, while {the} other conditions are satisfied {in Theorem \ref{(3.6)}}. {Now,} the {graph} of {$G^{-1}_{Y_{1:3}}\big(F_{X_{1:3}}(x)\big)/x$ is} shown in {Figure} \ref{(12)}(a), which is non-negative in $y$. {From this figure, it is clear that $X_{1:3}\nleq_{*}Y_{1:3}$.}
\end{counterexmp}

{The following example provides an illustration of the result in Theorem \ref{(3.7)}.}
\begin{exmp}\label{example5.7}
Set $n=3$, $\beta=0.5$, $(\alpha_1,\alpha_2,\alpha_3)=(0.1,0.2,0.3)$, and $(\alpha^*_1,\alpha^*_2,\alpha^*_3)=(0.4,0.6,0.8)$. Figure \ref{(12)}(b) represents the curve of {$G^{-1}_{Y_{1:3}}\big(F_{X_{1:3}}(x)\big)$} is positive and this means that $X_{1:3}\leqslant_{c}Y_{1:3}$. Hence the result in Theorem \ref{(3.7)} can hold if $\alpha^*_k>\alpha_k$.
\end{exmp}

\begin{figure}[H]
\begin{center}
{\includegraphics[width=3.8in,height=2.2in]{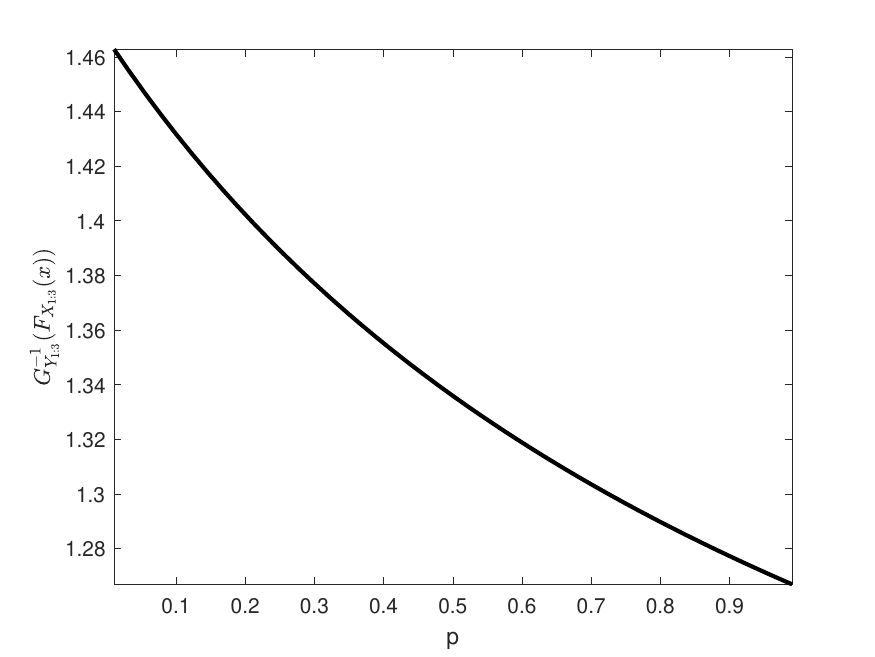}}
\caption{{Plot of $G^{-1}_{Y_{1:3}}\big(F_{X_{1:3}}(x)\big)$ as in Counterexample \ref{counexam5.7}.}}
\label{(14)}
\end{center}
\end{figure}

{We now present a counterexample to emphasize that the condition ``$\alpha^*_k\geq\alpha_k$'' is required for the result in Theorem \ref{(3.7)}.}
\begin{counterexmp}\label{counexam5.7}
Set $n=3$, $\beta=2$, $(\alpha_1,\alpha_2,\alpha_3)=(4.4,4.6,5.8)$, and $(\alpha^*_1,\alpha^*_2,\alpha^*_3)=(3.1,3.2,3.3)$. Here, $(4.4,4.6,5.8)\nleq(3.1,3.2,3.3)$, while {the} other conditions are satisfied {in Theorem \ref{(3.7)}}. In Figure \ref{(14)}, the {graph} of $G^{-1}_{Y_{1:3}}\big(F_{X_{1:3}}(x)\big)$ is not negative and this means that $X_{1:3}\nleq_{c}Y_{1:3}$. Hence the result in Theorem \ref{(3.7)} does not hold for $\alpha_k\nleq\alpha^*_k$.
\end{counterexmp}

\section{Conclusion}
Stochastic comparisons of heterogeneous series systems when component follows $\mathcal{LFR}$ distributed components under various parameter settings were studied for using usual stochastic, hazard rate, and likelihood ratio orderings. Also, investigated stochastic comparison results for the usual stochastic order for parallel systems. For series systems, stochastic comparison results have been studied for transform and variability orders such as dispersive, convex, star, Lorenz orders. All results have been illustrated with examples and counterexamples using graphics. The established results will be useful in the selection of a more reliable system. Stochastic comparison findings have been examined for independent $\mathcal{LFR}$ distributed component lifetimes. This problem can also be explored for the dependent heterogeneous components. In addition, consider a real-life example of a water distribution system in a residential area where water flows through multiple pipelines. In a series configuration, water must pass sequentially through three pipes $P_1$, $P_2$, and $P_3$, each with different aging characteristics, represented by heterogeneous linear failure rates: $\lambda_1(t)=0.01t$, $\lambda_2(t)=0.02t$, and $\lambda_3(t)=0.03t$, where $t$ is the time in years. The reliability of the series system is the product of the individual reliabilities:
\begin{eqnarray*}
R_s(t)=e^{-\int_{0}^{t}(0.01u+0.02u+0.03u)du}=e^{-0.03t^2}.
\end{eqnarray*}
In contrast, if the system is designed in a parallel configuration, water can flow through any of the three pipes independently, ensuring supply as long as at least one pipe is functional. The reliability of the parallel system is given by:
\begin{eqnarray*}
R_p(t)=1-\prod_{i=1}^{3}(1-e^{-\int_{0}^{t}\lambda_i(u)du})=1-(1-e^{-0.005t^2})(1-e^{-0.01t^2})(1-e^{-0.015t^2}).
\end{eqnarray*}
Stochastically, the parallel system is more reliable since it can withstand individual pipe failures, whereas the series system fails if any single pipe fails, leading to a much faster decay in reliability over time.	
\bibliographystyle{unsrt}
\bibliography{ref}

\begin{thebibliography}{10}

\bibitem{kuo2000annotated}
Way Kuo and V~Rajendra Prasad.
\newblock An annotated overview of system-reliability optimization.
\newblock {\em IEEE Transactions on Reliability}, 49(2):176--187, 2000.

\bibitem{balakrishnan19981}
N~Balakrishnan and Calyampudi~R Rao.
\newblock Order {S}tatistics: An {I}ntroduction.
\newblock {\em Handbook of Statistics}, 16:3--24, 1998.

\bibitem{david2004order}
Herbert~Aron David and Haikady~Navada Nagaraja.
\newblock Order {S}tatistics.
\newblock {\em Encyclopedia of Statistical Sciences}, 2004.

\bibitem{shaked2007stochastic}
Moshe Shaked and J~George Shanthikumar.
\newblock {\em Stochastic Orders}.
\newblock Springer Science \& Business Media, 2007.

\bibitem{jewitt1991}
I~Jewitt.
\newblock Applications of likelihood ratio orderings in economics.
\newblock {\em Lecture Notes-Monograph Series}, pages 174--189, 1991.

\bibitem{nanda2001hazard}
Asok~K Nanda and Moshe Shaked.
\newblock The hazard rate and the reversed hazard rate orders, with
  applications to order statistics.
\newblock {\em Annals of the Institute of Statistical Mathematics},
  53(4):853--864, 2001.

\bibitem{balakrishnan2014stochastic}
Narayanaswamy Balakrishnan, Abedin Haidari, and Khaled Masoumifard.
\newblock Stochastic comparisons of series and parallel systems with
  generalized exponential components.
\newblock {\em IEEE Transactions on Reliability}, 64(1):333--348, 2014.

\bibitem{gupta2015stochastic}
Nitin Gupta, Lakshmi~Kanta Patra, and Somesh Kumar.
\newblock Stochastic comparisons in systems with {F}r{\`e}chet distributed
  components.
\newblock {\em Operations Research Letters}, 43(6):612--615, 2015.

\bibitem{fang2015stochastic}
Longxiang Fang and Xinsheng Zhang.
\newblock Stochastic comparisons of parallel systems with exponentiated
  {W}eibull components.
\newblock {\em Statistics $\&$ Probability Letters}, 97:25--31, 2015.

\bibitem{chowdhury2017stochastic}
Shovan Chowdhury and Amarjit Kundu.
\newblock Stochastic comparison of parallel systems with log-{L}indley
  distributed components.
\newblock {\em Operations Research Letters}, 45(3):199--205, 2017.

\bibitem{fang2012new}
Longxiang Fang and Xinsheng Zhang.
\newblock New results on stochastic comparison of order statistics from
  heterogeneous {W}eibull populations.
\newblock {\em Journal of the Korean Statistical Society}, 41(1):13--16, 2012.

\bibitem{misra2012new}
Neeraj Misra and Amit~Kumar Misra.
\newblock New results on stochastic comparisons of two-component series and
  parallel systems.
\newblock {\em Statistics $\&$ Probability Letters}, 82(2):283--290, 2012.

\bibitem{torrado2015stochastic}
Nuria Torrado and Subhash~C Kochar.
\newblock Stochastic order relations among parallel systems from {W}eibull
  distributions.
\newblock {\em Journal of Applied Probability}, 52(1):102--116, 2015.

\bibitem{kundu2016ordering}
Amarjit Kundu and Shovan Chowdhury.
\newblock Ordering properties of order statistics from heterogeneous
  exponentiated {W}eibull models.
\newblock {\em Statistics $\&$ Probability Letters}, 114:119--127, 2016.

\bibitem{balakrishnan2018ordering}
Narayanaswamy Balakrishnan, Phalguni Nanda, and Suchandan Kayal.
\newblock Ordering of series and parallel systems comprising heterogeneous
  generalized modified {W}eibull components.
\newblock {\em Applied Stochastic Models in Business and Industry},
  34(6):816--834, 2018.

\bibitem{patra2018some}
Lakshmi~Kanta Patra, Suchandan Kayal, and Phalguni Nanda.
\newblock Some stochastic comparison results for series and parallel systems
  with heterogeneous {P}areto type components.
\newblock {\em Applications of Mathematics}, 63(1):55--77, 2018.

\bibitem{kayal2022stochastic}
Suchandan Kayal and Phalguni Nanda.
\newblock Stochastic comparisons of parallel systems with generalized
  {K}umaraswamy-{G} components.
\newblock {\em Communications in Statistics-Theory and Methods},
  51(14):4712--4738, 2022.

\bibitem{kayal2023some}
Suchandan Kayal, Raju Bhakta, and N~Balakrishnan.
\newblock Some results on stochastic comparisons of two finite mixture models
  with general components.
\newblock {\em Stochastic Models}, 39(2):363--382, 2023.

\bibitem{barmalzan2022usual}
Ghobad Barmalzan, Sajad Kosari, and Narayanaswamy Balakrishnan.
\newblock Usual stochastic and reversed hazard orders of parallel systems with
  independent heterogeneous components.
\newblock {\em Communications in Statistics-Theory and Methods},
  51(14):4781--4806, 2022.

\bibitem{kodlin1967new}
D~Kodlin.
\newblock A new response time distribution.
\newblock {\em Biometrics}, pages 227--239, 1967.

\bibitem{bain1974analysis}
Lee~J Bain.
\newblock Analysis for the linear failure-rate life-testing distribution.
\newblock {\em Technometrics}, 16(4):551--559, 1974.

\bibitem{sen1995inference}
Ananda Sen and Gouri~K Bhattacharyya.
\newblock Inference procedures for the linear failure rate model.
\newblock {\em Journal of Statistical Planning and Inference}, 46(1):59--76,
  1995.

\bibitem{muller2002comparison}
Alfred M{\"u}ller and Dietrich Stoyan.
\newblock Comparison methods for stochastic models and risks.
\newblock 2002.

\bibitem{li2013stochastic}
Haijun Li and Xiaohu Li.
\newblock Stochastic orders in reliability and risk.
\newblock {\em Honor of Professor Moshe Shaked. Springer, New York}, 2013.

\bibitem{kleiber2002variability}
Christian Kleiber.
\newblock Variability ordering of heavy-tailed distributions with applications
  to order statistics.
\newblock {\em Statistics $\&$ {P}robability {L}etters}, 58(4):381--388, 2002.

\end{thebibliography}
\end{document}